\title {      Values of Brownian   
             intersection exponents I: \\ 	
           Half-plane exponents  
}
\author { Gregory F.  Lawler\thanks{Duke University}
\and Oded Schramm\thanks{The Weizmann Institute of Science and
Microsoft Research}
\and Wendelin Werner\thanks{Universit\'e Paris-Sud}
}
\date {}

\documentclass[12pt]{article}
\usepackage{amssymb}
\usepackage{amsmath}
\usepackage{amsthm}
\usepackage{amsfonts}

\numberwithin{equation}{section}
\newtheorem{theorem}{Theorem} 
\newtheorem{corollary}[theorem]{Corollary} 
\numberwithin{theorem}{section} 
\newtheorem{lemma}[theorem]{Lemma} 
\newtheorem{prop}[theorem]{Proposition}

\newcommand{\R}{\mathbb{R}}
\newcommand{\C}{\mathbb{C}}
\newcommand{\Z}{\mathbb{Z}}
\newcommand{\N}{\mathbb{N}}
\newcommand{\HH}{\mathbb{H}}
\def \uu {{s}}
\def \eps {\epsilon}
\def \tx {{\tilde \xi}}
\def \P {{\bf P}}
\def \prob {{\bf P}}
\def \expect {{\bf E}}
\def \p {{\partial}}

\def \proof {{ \medbreak \noindent {\bf Proof.} }}

\def\Ppb#1{\P \bigl[#1\bigr]}

\def \st {:}
\def\slepar{\kappa}
\def\SG/{$SLE_\slepar$}
\def\Wsle{W^\slepar}
\def\Ss/{$SLE_6$}

\def\closure#1{{\overline{#1}}}
\def\diam{{\rm diam}} 
\def\defeq{:=}
\def\term#1{{\bf #1}}
\def\Im{{\rm Im}} 
 
\def\Ua#1{U^2(#1)}
\def\Ub#1{U^1(#1)}

\begin {document}

\maketitle

\begin {abstract}
This paper proves conjectures originating in the
physics literature regarding the intersection exponents
of Brownian motion in a half-plane.  
For instance, suppose that $B$ and $B'$ are two
independent planar Brownian motions started from 
distinct points in a half-plane ${\cal H}$. Then 
as $t \to \infty$,
$$
P \Bigl[ B[0,t] \cap B' [0,t] = \emptyset \hbox { and }
B[0,t] \cup B' [0,t] \subset  {\cal H} \Bigr]
= t^{-5/3 + o(1)}.
$$
The proofs use ideas and tools 
developed by the authors in 
previous papers. We prove that one of the 
stochastic L\"owner evolution 
processes (with parameter 6, that we will call
\Ss/ and which has been conjectured to 
correspond to the scaling limit of critical percolation cluster 
boundaries) 
satisfies  the ``conformal restriction property''.  
We establish a  generalization of 
Cardy's formula (for crossings of a rectangle 
by a percolation cluster) for \Ss/, 
from which the exact values of
intersection exponents for \Ss/  follow.
Since this process satisfies the conformal restriction
property, 
the Brownian intersection exponents can be determined from the 
\Ss/ intersection exponents.

Results about intersection
exponents in the whole plane will appear in subsequent papers.
\end {abstract}

\section {Introduction}

Theoretical physics predicts that
conformal invariance plays a crucial role
in the macroscopic behavior
of  a wide
class of two-dimensional models in statistical physics
(see, e.g., \cite {BPZ, Ca1}).
For instance, by making the assumption 
that critical planar percolation behaves in a 
conformally invariant way in the scaling limit
and using ideas involving  conformal field theory,
 Cardy \cite {CaFormula}
produced an exact formula for the limit, as $N\to\infty$,
of the probability 
that, in two-dimensional critical percolation,
there exists a cluster crossing   
the rectangle $[0, aN] \times [0, b N]$.
Also, Duplantier and Saleur \cite {DS2} 
predicted the ``fractal dimension'' of the hull of a very large
percolation cluster.
These are just two examples among many such predictions.

In 1988, Duplantier and Kwon \cite {DK} suggested that the 
ideas of conformal field theory can also be applied to predict
the intersection exponents between random walks in $\Z^2$
(and Brownian motions in $\R^2$).
They predicted, for instance, that if $B$ and $B'$ are
independent 
planar Brownian motions
(or simple random walks in $\Z^2$)  started  
 from distinct points
in the upper half-plane
 $\HH = \{ (x,y) \st 
y >  0 \}= \{z\in\C\st \Im(z)>0\}$, then 
when $n \to \infty$,
\begin {equation}
\label {01}
\Ppb { B[0, n] \cap B' [0,n] = \emptyset } 
=
 n^{-\zeta + o (1)} 
\end {equation}
and
\begin {equation}
\label {02}
\Ppb  
{ B[0,n] \cap B' [0, n ] = \emptyset 
\hbox { and }
B [0,n] \cup B'[0, n] \subset \HH  
}
=
n^{-\tilde \zeta  + o (1)}
\end  {equation}
where
$$ \zeta = 5/8\,,\qquad\tilde \zeta = 5/3. 
$$
Very recently, Duplantier \cite {Dqg} 
gave another physical derivation of these conjectures based 
on ``quantum gravity''. 

In 1982, Mandelbrot \cite {M}  suggested that the 
Hausdorff dimension of the Brownian frontier (i.e., the boundary of
a connected component of the complement
of the path) is 4/3, based on simulations and the analogy
 with the conjectured
value for the fractal dimension of self-avoiding walks 
predicted by Nienhuis (also 4/3; see, e.g., \cite {MS}).

\medbreak

To date, none of the physicists' arguments have been made rigorous,
and it seems very difficult to use their methods
to produce proofs.
Very recently, Kenyon \cite {K1,K2,K3}  
managed to derive the 
exact values of  critical exponents for  ``loop-erased
random walk'' that 
theoretical physicists had predicted
(Majumdar \cite {Maj},
 Duplantier \cite {Dle}). Kenyon's  
 methods
 involve the relation of the loop-erased walk to the uniform
 spanning tree and to
domino tilings. 
Kenyon shows that the equations relating probabilities
of some domino tiling events are discrete analogues
of the Cauchy-Riemann equations, and therefore the probabilities
can be approximated by analytic functions with prescribed boundary
behavior.  These methods do not seem applicable for the goals
of the present paper.

For planar Brownian motions, it is
 easy to show, using subadditivity arguments 
and the scaling property, that there exist positive finite 
numbers $\zeta$ and $\tilde \zeta$ 
such that 
(\ref {01}) and (\ref {02}) are true.
Up to the present paper, there was not even a mathematical
heuristic arguing that
the values of $\zeta$ and $\tilde \zeta$ are  $5/8$ and $ 5/3$.
Burdzy-Lawler \cite {BL} (see also \cite {CM, LP}) 
showed that the intersection exponents were indeed the same
for simple random walks as for Brownian motions; 
Lawler~\cite {Lcut}  proved that
the Hausdorff dimension of  the
set of cut points of a 
Brownian path  is $2 - 2 \zeta$.
He also showed
 (see \cite {Lfront, Lmulti}) 
that the Brownian frontier
(and more generally the whole multifractal spectrum of the Brownian
frontier) can be expressed in terms of 
exponents defined analogously to 
$\zeta$.  As part of that 
 work, he showed that 
the right hand side of (\ref{01}) can be replaced with 
$n^{-\zeta} g(n)$ where $g$ is bounded away from $0$ and
infinity; we expect that the argument
can be adapted to show that 
 the same is true
for (\ref{02}).

Recently, Lawler and Werner \cite {LW1} extended
the definition of intersection exponents in a natural
way
to ``non-integer packets of Brownian motions'' 
and  derived certain functional relations between these exponents.
These relations indicate 
 that Mandelbrot's conjecture  that the dimension of the
Brownian frontier is $4/3$ 
is indeed compatible with the predictions of Duplantier-Kwon.
It turned out that 
intersection
exponents in the half-plane play an important  role 
in understanding exponents in the whole plane.
Conformal invariance of planar Brownian motion is a crucial
tool in the derivation of these relations.
In particular, there is a measure on
Brownian excursions in domains that has
some strong conformal invariance properties,
including a ``restriction'' (or ``locality'')
property.

In another paper, Lawler and Werner \cite {LW2} showed  
that intersection exponents associated to
 any  conformally invariant measure on sets  
with this restriction
 property are very closely related to the 
Brownian exponents.
This provides a rigorous justification 
to the link between the conjectures regarding 
intersection exponents for planar Brownian motions 
 and conjectures for 
intersection exponents 
of critical 
 percolation clusters (see \cite {DS2, Ca2, ADA}), 
because  percolation clusters 
are conjectured to be conformally 
invariant in the scaling limit --- see, e.g., \cite {LPS,Ai} 
 ---  and they should
also  have a restriction property (because of the 
independence properties of percolation). 
The question
of how to compute these exponents  
remained open.

Independently, Schramm \cite {S1} defined a new class of 
conformally invariant stochastic processes indexed
by a real parameter $\slepar \ge 0$, called \SG/
(for stochastic L\"owner evolution process  
with parameter
$\slepar$). The definition of these processes 
is based on 
L\"owner's ordinary differential equation that encodes
in a conformally invariant way a continuous family of
shrinking domains (see, e.g., \cite {Lo, P1}).
More precisely, \cite{S1} defines a family of 
conformal maps $g_t$ from  subsets $D_t$ 
of $\HH$ onto $\HH$
by the equation 
\begin{equation}\label{introlow} 
\partial_t g_t (z) 
=  
\frac {-2}{\beta_{\slepar t} - g_t (z)}\,,
\end{equation}
where $\beta$ is a standard
Brownian motion on the real line.
(Actually, in~\cite{S1}, instead of~(\ref{introlow}),
the corresponding equation for the inverse maps $g_t^{-1}$
is considered.) 
The domain $D_t$
can be defined as the set of $z_0\in\HH$ such that 
a solution $g_{s}(z_0)$ of this equation exists for
$s\in[0,t]$.
When $t$ increases, the set 
$K_t= \HH \setminus D_t$ increases: Loosely
speaking, $(K_t, t \ge 0)$ can be viewed as a 
 growing ``hull'' that is 
penetrating the half-plane.
By applying a conformal homeomorphism $f:\HH\to D$,
\SG/ can similarly be defined
in any simply connected domain $D\subsetneqq\C$.

In \cite {S1}, the main focus is
on the case $\slepar=2$, which is conjectured there to correspond
to  the scaling limit of loop-erased random walks, 
but the conjecture that \Ss/
corresponds to the scaling limit of 
critical percolation cluster boundaries is also mentioned.
In particular (see \cite {S2}), 
it is 
possible  to compute explicitly the probability that an
\Ss/ crosses a rectangle of size $a \times b$. It 
turns out that this result is exactly Cardy's formula.
This gives a mathematical proof for Cardy's formula,
assuming the still open conjecture that
\Ss/ is indeed the scaling limit of percolation cluster 
boundaries. 

\medbreak

The main goal of the present paper is to prove some of the 
conjectured values of intersection exponents
of Brownian motion in a half plane.  

\begin {theorem}
\label {theoremintro}
Let $B^1, \cdots, B^p$ denote $p$ independent planar 
Brownian motions ($p \ge 2$) started from distinct points in the 
upper half-plane $\HH$. Then, when $t \to \infty$,
\begin {equation*}
\Ppb
{ \forall  i \not= j \in \{1, \cdots , p\},  \
B^i [0,t] \cap B^j [0,t] = \emptyset
\hbox { and } B^i [0,t] \subset \HH 
}
= t^{- \tilde \zeta_p + o (1) }\,,
\end {equation*}
where
$$
\tilde \zeta_p = \frac {p (2p+1)}{6}
\,.$$
\end  {theorem}

These values have been predicted by Duplantier and Kwon \cite {DK}.
In particular $\tilde \zeta=  \tilde{\zeta_2 }= 5/3$.

We also establish the exact value (and confirm  some of
the conjectures stated in \cite {LW1,Dqg}) 
of more general intersection exponents between packets of 
Brownian motions in the  half-plane; see Theorem~\ref {general}.

\medskip

The proof of  Theorem~\ref{theoremintro} uses a combination of ideas
from the papers \cite {LW1, LW2, S1, S2}.
However, to make the paper more accessible and self-contained, we attempt to
review and explain all the necessary background.  The reader who wishes
to see complete proofs for all stated theorems has to be familiar
with the basics of stochastic calculus and conformal mapping theory,
and read about the excursion measure and the cascade relations from
\cite{LW1}.

\medbreak

Although, at present, a proof of the conjecture that \Ss/ is the scaling
limit
of critical percolation cluster boundaries seems out of reach,
this conjecture does lead
one to believe that \Ss/ must satisfy a ``locality'' property, namely,
it is not affected by the boundary of a domain when it is
in the interior.  This locality property for \Ss/ is stated
more precisely and proved in Section~\ref{s.loc}.
It is worthwhile to note that the locality property does not
hold for the \SG/ processes when $\slepar\neq6$.

In Section~\ref{s.Sexp}, we 
prove that \Ss/ satisfies  a generalization of Cardy's formula
for percolation crossings probabilities.  From this, exponents associated
with
the \Ss/ process are computed.  

In Section~\ref{s.Bexp}, universality ideas from~\cite{LW2} are used to
compute
the half-plane Brownian exponents from the \Ss/ exponents, which 
completes the proof of Theorem~\ref{theoremintro}.

In a final short Section~\ref{s.percolation}, the conjectured 
relationship between \Ss/ and critical percolation is discussed. 
It is demonstrated that this conjecture implies a formula
from the physics literature \cite{DS2, Ca2, ADA} for the exponents
corresponding
to the event that there are $k$ disjoint
percolation crossings of a long rectangle.

\medskip

In the subsequent
papers \cite {LSW2,LSW3,LSWan}, we determine the exponents 
in the full plane and the remaining half-plane exponents.
In particular, we prove that $\zeta=5/8$,
and also establish Mandelbrot's conjecture that the Hausdorff dimension
of the frontier of planar Brownian motion is $4/3$.

\medskip

It might be worthwhile to explain why
 the Brownian intersection exponents
 are accessible through \Ss/, but are difficult 
to compute directly.  In a way, the \Ss/
process
is simpler, since $K_t$ continuously grows from its outer boundray.  This
means that when studying its evolution, one can essentially forget
its interior, and only keep track of the exterior of $K_t$.  By conformal
invariance, this reduces problems to finitely many dimensions.  The
situation
with planar Brownian motion is completely different, since it may enter
holes
it has surrounded and emerge to the exterior someplace else.
Many computations with \SG/ are readily convertible to PDE problems,
and in the presence of enough symmetry, some variables can often be
eliminated, converting the PDE to an ODE.

\section {\Ss/ and its locality property}
\label {s.loc}

\subsection {The definition of chordal \SG/ and some basic properties}
\label{slebasic}

Let $(\beta_t, t \ge 0)$ be a standard
real-valued Brownian motion starting at $\beta_0=0$,
let $\slepar>0$, and let $\Wsle_t=\beta_{\slepar t}$.
Consider the ordinary differential equation
\begin{equation}\label{e.low}
\p _t g_t(z) = \frac {-2}{\Wsle_t - g_t(z)}
\end{equation}
with $g_0(z)=z$.
For every $z_0\in\HH$ and every $T>0$,
either there is a solution of (\ref{e.low})
for $t\in[0,T]$ and for all $z$ in a neighborhood
of $z_0$, or there is some $t_0\in(0,T]$ such that
the solution exists for $t\in[0,t_0)$ and
$\lim_{t\nearrow t_0} g_t(z)= \Wsle_{t_0}$.
Let $D_T$ be the (open) set of $z\in\HH$ such
that the former is true, and let $K_T$ be the
set of $z\in\HH$ such that the latter holds.
By considering the inverse flow
$\p _t G_t(z)= 2\bigl(\Wsle_{T-t}-G_t(z)\bigr)^{-1}$,
it is easy to see that $g_t(D_t)=\HH$, and that 
$g_t:D_t\to\HH$ is conformal.
The process $g_t$, $t\geq 0$, will be called
the \term{chordal stochastic L\"owner evolution process}
with parameter $\slepar$, or just \SG/;
see \cite{S1,S2}.
In \cite{S1}, a variation of this process, which we now call
\term{radial \SG/} was also studied.  In the current paper, we will not
use radial \SG/, and therefore the word ``chordal'' will usually be 
omitted.  (However, radial \SG/ plays a major role in a
subsequent paper \cite{LSW2}.)
The set $K_t=\HH\setminus D_t$ will be called the \term{hull} of the SLE.
The process $\Wsle_t$ will be called the \term{driving process}
of the SLE.

It is easy to verify that each of the maps $g_t$
satisfies the \term{hydrodynamic normalization}
at infinity:
\begin{equation}\label{hydro}
\lim_{z\to\infty}g(z)-z = 0\,.
\end{equation}

\bigskip
\noindent{\bf Remarks.}
It will be shown \cite{RS} that for all $\kappa\ge 0$
the hull $K_t$ of \SG/ is a.s.\ generated by a path.
More precisely, a.s.\
the map $\gamma(t)\defeq g_t^{-1}(\Wsle_t)$ is a well
defined continuous path in $\overline\HH$
and for every $t\ge 0$ the domain $D_t$
is the unbounded connected component of $\HH\setminus\gamma([0,t])$.
There, it will also be shown that when $\slepar\le4$
a.s.\ $K_t$ is a simple path for all $t>0$.
  This is not the case when $\slepar>4$ \cite{S1}.
However, these results will not be needed for 
the current paper or for \cite{LSW2,LSW3,LSWan}.

L\"owner \cite{Lo} considered the equation
$$
\p _t g_t(z) = g_t(z)\frac{\zeta(t)+g_t(z)}{\zeta(t)-g_t(z)}
\,,
$$
with $g_0(z)=z$, where $z$ is in the unit disk, and $\zeta(t)$
is a parameter.  He used this equation in the
study of extremal problems
for classes of normalized conformal mappings.
In L\"owner's differential equation, the maps $g_t$
satisfies  $g_t(0) = 0$.  The equation (\ref{e.low}) is an analogue
of L\"owner's equation in the half plane, where the
boundary point $\infty$ is fixed instead of  $0$ and $\zeta(t)$
is chosen to be scaled Brownian motion.

Marshall and Rohde \cite{MR} study conditions on $\zeta(t)$ which
imply that $K_t$ is a simple path.
\bigskip

We now note some basic properties of \SG/.

\begin{prop}

(i)
[Scaling] \label {scaling}
\SG/ is scale-invariant in the following
sense.  Let $K_t$  be the hull of \SG/, and let $\alpha>0$.
Then the process $t\mapsto\alpha^{-1/2} K_{\alpha t}$
has the same law as $t\mapsto K_t$.

(ii)
[Stationarity] \label{p.stationary}
Let $g_t$ be an \SG/ process in $\HH$, driven
by $\Wsle_t$, and let $\tau$ be any stopping time. 
Set $\tilde g_t(z)=g_{\tau+t}\circ g_{\tau}^{-1}(z+\Wsle_\tau)-\Wsle_\tau$.
Then $\tilde g_t$ is an \SG/ process in $\HH$ starting
at $0$, which is independent from
$\{g_t\st t\in[0,\tau]\}$.
\end{prop}

\proof

(i) If $K_t$ is driven by $\Wsle_t$, then 
$\alpha^{-1/2} K_{\alpha t}$ is
driven by $\alpha^{-1/2} \Wsle_{\alpha t}$,
which has the same law as $\Wsle_t$.

(ii) The process 
$\tilde g_t$ is driven by $\Wsle_{t+\tau}-\Wsle_\tau$.
\qed\bigskip

\bigskip
We now consider the definition of \SG/ in domains other than $\HH$.
\medskip

Let $f:D\to\HH$ be a conformal homeomorphism from
some simply connected domain $D$.
Let $f_t$ be the solution of (\ref{e.low}) with $f_0(z)=f(z)$.
Then $(f_t, t \ge 0)$ will be called the \SG/ in $D$ starting at $f$.
If $g_t$ is the solution of (\ref{e.low}) with $g_0(z)=z$,
then we have $f_t=g_t\circ f$. If $K_t$ is the hull associated to
$g_t$, then  the hull  associated
with $f_t$ is just $f^{-1}(K_t)$.

Suppose that $\partial D$ is a Jordan curve in $\C$,
and let $a,b\in\partial D$ be distinct.
Then we may find such an 
$f:D\to\HH$ with $f(a)=0$ and $f(b)=\infty$. 
Let $K^f_t$ be the \SG/ hull associated with the \SG/ process starting
at $f$.  If $f^*$ is another such map $f^*:D\to\HH$ with
$f^*(a)=0$ and $f^*(b)=\infty$, then
$f^*(z) = \alpha f(z)$ for some $\alpha > 0$.
By Proposition  \ref{scaling},
the
corresponding \SG/ hull $K^{f^*}_t$ has the same law as a linear
time-change of $K^f_t$.
This makes it natural to consider $K^f_t$ as a process from
$a$ to $b$ in $D$, and to ignore the role of $f$. 
However, when $D$ is not a Jordan curve, some care may
be needed since the conformal map $f$ does not necessarily extend
continuously to the boundary.  Partly for that reason, we have
chosen to stress the importance of the conformal parameterization $f$.

\subsection{The locality  property}

The main result of this section can be loosely 
described as follows: an \Ss/ process  does not feel
where the boundary of the domain lies as long as it 
does not hit it. This is consistent
with the conjecture \cite {S1}
that the \Ss/ process is  
the scaling limit of percolation cluster boundaries,
which is explained in Section~\ref{s.percolation}.  
This restriction property  can therefore
be viewed as additional evidence in favor of  this  conjecture.
This feature is special to \Ss/; it is
not shared by \SG/ when $\slepar\neq 6$.

Such properties were studied in \cite {LW2} and 
called ``complete conformal invariance''  
(when combined with a conformal invariance property).
As pointed out there, all processes 
with complete conformal invariance 
have closely related intersection exponents.

\def \nice {nice }  

Let us first state a general local version of 
this result.
We will say that the path $\gamma$ is \nice if it is a continuous 
simple  path $\gamma: [0,1] \to \overline \HH$, such that 
$\gamma(0), \gamma(1) \in \R \setminus \{0\} $
and $\gamma (0,1) \subset \HH$.
We then call the connected component $N= N(\gamma)$ of $\HH
\setminus \gamma [0,1]$ such that $0 \in \p N$
a \nice neighborhood of $0$ in $\HH$.
Note that $N$ can be bounded
or unbounded, depending on the sign of $\gamma(0) \gamma (1)$.
When $N$ is a \nice neighborhood of $0$, one can  find 
a conformal homeomorphism $\psi= \psi_N$ from $N$ onto $\HH$
such that $\psi (0)= 0$, $\psi' (0) = 1$
 and $\psi^{-1} (\infty)$ is equal to $\infty$ if $N$ is unbounded and
to $\gamma (1)$ if $N$ is bounded.

\begin {theorem}[Locality]  
 \label{t.loc}
Let $f:D\to\HH$ be a conformal homeomorphism
from a domain $D\subset\C$ onto $\HH$. 
Suppose that $N$ is a \nice neighborhood of $0$
in $\HH$.
Define 
$D^* = f^{-1} (N)$ 
and let $f^*$ be the conformal homeomorphism $\psi_N \circ f $
from $D^*$ onto $\HH$.
Let $K_t\subset D$ be the hull of \Ss/
starting at $f$, and let
$\tau\defeq \sup \{t\st \closure{K_t}\cap \p D^*\cap D=\emptyset\}$. 
Let $K_t^*$ denote \Ss/  in $D^*$ started at 
$f^*$ 
and let $\tau^*\defeq\sup\{t\st \closure{K^*_t}\cap \p D^*\cap
D=\emptyset\}$.
 
Then 
the law of $(K_t, t < \tau)$ is that of a time-change of 
$(K_t^*, t < \tau^*)$.
\end {theorem}

Note that in this theorem, we have not made any regularity assumption 
on the boundary of the domain $D$. 

A consequence of this result is that, modulo time-change, one can define
the hull of \Ss/ in a non-simply connected domain with
finitely many boundary components since such   
a domain looks locally like a simply connected domain.

\medbreak

This property implies the following 
``global'' 
restriction properties.
For convenience only, we will state them
under some assumptions  on the boundaries of the domains.

\begin {corollary} [Splitting property]
\label {split}
Let $D$ denote a simply connected domain such that 
$\p D$ is a Jordan curve. Let $a$, $b$ and $b'$ denote 
three distinct points on $\p D$, and let  $I$ 
denote the connected component of $\p D \setminus \{ b,b'\}$
that does not contain $a$.
Let $(K_t, t \ge 0)$
 (respectively $K_t'$) denote an \Ss/ in $D$ 
from $a$ to $b$ (resp. from $a$ to $b'$). Let $T$ (resp. $T'$)
denote the first time at which $\closure{ K_t}$ (resp. $\closure{ K_t'}$)
intersects $I$.
Then $( K_t , t < T )$ and $ ( K_t'  , t < T' )$ 
have the same law up to time-change.
\end {corollary}
 
\begin {corollary} [Restriction property]
\label {rest}   
Let $D^* \subset D  $ denote two simply connected domains,
and assume that $\p D$ is a Jordan curve.
Suppose that
$I \defeq \p D^* \setminus \p D$
is connected.
Take two distinct points  
$a$ and $b$ in $ \p D \cap \p D^* \setminus \closure I$.

Let $(K_t, t \ge 0)$ denote \Ss/ from $a$ to $b$ in $D$, and 
$T\defeq \sup\{t\st\closure{K_t}\cap I=\emptyset\}$. 
Similarly, let $(K_t^*, t \ge 0)$ be \Ss/
from $a$ to $b$ in $D^*$, and 
$T^*\defeq \sup\{t\st\closure{K_t^*}\cap I=\emptyset\}$. 
Then,
$( K_t , t < T )$ and $( K_t^* , t < T^* )$
have the same law up to time-change.
\end {corollary}

In the present paper, we will use these results when $D$
is a rectangle.
\medbreak

\noindent
{\bf Proof of Corollary \ref {split}} (assuming Theorem \ref {t.loc}).
This is just a consequence of the fact that in Theorem \ref {t.loc} with
$D= \HH$ and bounded $N$, one
can  replace $\gamma$ by $\beta (s) \defeq \gamma (1-s)$. Then we 
get that  the law of 
\Ss/ in $N$ from $0$ to $\gamma(0)$ is that of a time-change of \Ss/ in 
$N$ from $0$ to $\gamma (1)$ up to their hitting times
of $\gamma$. The result in a general domain follows by 
mapping it conformally onto  a nice neighborhood $N$ with $a$ mapped 
to $0$ and $\{b,b'\}$ to $\{\gamma(0),\gamma(1)\}$. 
\qed
 
\medbreak

\noindent 
{\bf Proof of Corollary \ref {rest}} (assuming Theorem \ref {t.loc}).
By approximation, it suffices to consider the
case where $I$ is a simple path.
Let $f$ denote a conformal map from $D$ onto $\HH$, 
with $f(a) = 0$ and $f(b) = \infty$.
Define $\gamma$ in such a way that 
 $\gamma [0,1] = f(I)$; note that $D^* = f^{-1} (N(\gamma))$.
As $b \in \p D^* \setminus I$, $N(\gamma)$ is unbounded.
Hence, by Theorem \ref {t.loc}, the law of \Ss/ in $D$
from $a$ to $b$ stopped when it hits $I$, is (up to time-change)
the same as that of \Ss/ in $D^*$ from $a$ to $b$ stopped when its
closure hits $I$.
\qed

\medbreak

In order to prove Theorem \ref {t.loc}, we will establish the following 
lemma:

\begin {lemma}
\label {l.loc}
Under the assumptions of Theorem \ref {t.loc},
define for any fixed $s<1$,
$L_s= \gamma (0,s]$, and 
$$ T= \sup \{ t \ge 0 \st
\overline K_t \cap \overline L_s = \emptyset \}.$$ 
For all $t \le T$,
let $g_{s,t}$ denote the conformal homeomorphism taking 
$\HH \setminus ( K_t \cup L_s) $ onto $\HH$
with 
the hydrodynamic normalization. 
Then, the process $(g_{s,t}, t < T)$ 
has the same law as a time-change of \Ss/ in $\HH\setminus L_s$ 
starting at $g_{0,s} - g_{0,s}(0)$, up to the time when the closure of its
hull
intersects $\closure{L_s}$. 
\end {lemma}

\noindent
{\bf Proof of Theorem \ref {t.loc}}
(assuming Lemma \ref {l.loc}).
In the setting of the Lemma, let 
$$ h_{s,t} (z) = \frac { 
g_{s,t} (z) - g_{s,0} (0)} { g_{s,0}' (0)}.
$$
By Lemma \ref {l.loc} and Proposition \ref {scaling},
$t \mapsto h_{s,t}$ has the same law as a time-change of \Ss/ 
starting at $h_{s,0}$.
Note that $h_{s,0} (0) = 0$, $h_{s,0}'(0) = 1$. 
Hence, it follows easily that for all $z \in N(\gamma)$,
$$
\lim_{s \to 1} h_{s,0} (z) = \psi_N (z) .
$$
By continuity, if we let $h_{1,t} = \lim_{s \to 1} h_{s,t}$,
then, $t \mapsto h_{1,t}$ 
has the same law as a time-changed  \Ss/ (in $N$) started from $\psi_N$.
The proof is completed by noting that the hull
of $h_{1,t}$ is $K_t$. 
\qed

\medbreak

The idea in the proof of Lemma \ref{l.loc} is to
study how the process $g_{s,t}$ changes as
$s$ increases.  For this, we will need to
use some of the properties of solutions to
(\ref{e.low}) where $\Wsle$ is replaced by
other continuous functions, and to study
how (deterministic)
families of conformal maps can be represented in this
way with some driving function.

\subsection{Deterministic expanding hulls}

\subsubsection {Definition and first properties}

If $(U_t, t \in [0, a])$ is a continuous  real-valued function,
then the process defined by 
\begin{equation}\label{genlow}
\p _t g_t(z) = \frac {-2}{U_t - g_t(z)}
\end{equation}
and $g_0(z)=z$
will be called the \term{L\"owner evolution with driving
function} $U_t$.
Note that $g_t$ satisfies the hydrodynamic normalization
(\ref{hydro}).
Moreover, 
\begin{equation}\label{expan}
g_t(z)=z+2tz^{-1}+a_2(t)z^{-2}+\cdots,\qquad z\to\infty,
\end{equation}
for some functions $a_j(t)$, $j=2,3,\dots$.
As above, we let $D_t\subset\HH$ denote the domain of $g_t$,
and let $K_t\defeq \HH\setminus D_t$.  $K_t$ will
be called the \term{expanding hull} of the process $g_t$.

We now address the question of which processes
$K_t$ can appear as the expanding hull driven by a continuous
function $U_t$.
We say that a bounded set $K\subset\HH$ is a
\term{hull} if $\HH\setminus K$ is open 
and simply connected.
The Riemann mapping theorem tells us that for
each hull $K$,  there is a unique
conformal homeomorphism $g_K:\HH\setminus K\to \HH$,
which satisfies the hydrodynamic normalization
(\ref{hydro}).  Let 
$$
A(K)=A(g_K)\defeq\frac 12\lim_{z\to\infty} z(g_K(z)-z);
$$
that is, $g(z)=z+2A(g)z^{-1}+\cdots$, near $\infty$.
Observe that $A(K)$ is real, because $g_K(x)$ is
real when $x\in\R $ and $|x|$ is sufficiently large.
Moreover, $A(K)\ge 0$, because $\Im\bigl(z-g_K(z)\bigr)$ 
is a harmonic function which vanishes at infinity 
and has nonnegative boundary values.  Note that
$$
A(g\circ h)=A(g)+A(h)
$$
if $g$ and $h$ satisfy the
hydrodynamic normalization.
It follows that $A(K)\le A(L)$ when $K\subset L$,
since $g_{L}=g_{g_K(L\setminus K)}\circ g_K$.

The quantity $A(g)$ is similar to capacity, and
plays an analogous role for the equation (\ref{e.low})
as capacity plays for L\"owner's equation.

\begin{theorem}\label{driven}
Let $(K_t, t\in[0,a])$ be an increasing family of hulls.
Then the following are equivalent.
\begin{enumerate}
\item
For all $t \in [0,a]$, $A(K_t)=t$,
and for each $\epsilon>0$ there is a $\delta>0$
such that for each $t\in[0,a-\delta]$ there is a bounded connected
set $S\subset\HH\setminus K_t$
with
$\diam(S)<\epsilon$  
and such that  $S$ disconnects  $K_{t+\delta}\setminus K_t$ 
from infinity in $\HH\setminus K_t$. 

\item
There is some continuous $U:[0,a]\to\R$,
such that ${K_t}$ is driven by $U_t$.
\end{enumerate}
\end{theorem}

In \cite{P1} a similar theorem is proven for
L\"owner's  differential equation in the disk.

Note that $\hat K_t$ may change discontinuously, in the
Hausdorff metric, as $t$ increases.  For example,
consider $\hat K_t\defeq \{\exp(is)\st 0<s\le t\}$ when
$t<\pi$ and $\hat K_{\pi}\defeq \{z\in\HH\st |z|\le 1\}$
and $\hat K_{t+\pi}\defeq K_{\pi}\cup (-1,-1+it]$, $t> 0$, say,
and let $K_t\defeq\hat K_{\phi(t)}$ where $\phi$ is chosen
to satisfy $A(K_{\phi(t)})=t$.

\begin{lemma}\label{cauchy}
Let $r>0$, $x_0\in\R$, and let $K$ be a hull contained in the disk
$\{z\st |z-x_0|<r\}$. Then
\begin{equation*}
\left|
 g_K^{-1}(z)-z+
\frac {2A(K)}{{z-x_0}}
\right|
 \leq \frac {C r A(K)}{|z-x_0|^{2}}
\end{equation*}
for all $z\in\HH$ with $|z-x_0|>Cr$, where $C>0$ 
is an absolute constant.
\end{lemma}

\medbreak
\noindent
{\bf Proof of Lemma \ref {cauchy}.}
For notational simplicity, we assume that $x_0=0$.
Clearly, this does not entail any loss of generality.
By approximation, we may assume that $K$ has smooth
boundary.  
Let $I\subset \R$ be
the smallest interval in $\R$ containing  
$\{g_K(x)\st x\in \p K\cap \HH\}$,
and let $f\defeq g_K^{-1}$.  Let $f_I$ be the restriction
of $f$ to $I$. 
Let $f^*$ denote the extension of
$f$ to $\C\setminus I$, by Schwarz reflection.  The Cauchy
formula gives
$$
2 \pi i f^*(w)
= \int_{|z|=R} \frac{f^*(z)}{z-w} \;dz 
+ \int_I \frac{{f_I(x)}-\overline{f_I(x)}}{x-w}\; dx\,,
$$
provided that $R>|w|$, $R> \max\{|x|\st x\in I\}$,
and $w\in\C\setminus I$.
Since $f^*(z) =z-2A(K)z^{-1}+\cdots$ near $\infty$,
\[  \lim_{R \rightarrow \infty} \int_{|z| = R}
      \frac{f^*(z)}{z-w} \; dz = \lim_{R \rightarrow
   \infty} \int_{|z| = R} \frac{z}{z-w} \; dz 
      =  2 \pi i w . \]
Consequently, we have
$$
f^*(w)-w
=  \frac{1}{\pi}
\int_I \frac{\Im\bigl({f_I(x)}\bigr)}{x-w} \;dx
\,.
$$
Multiplying by $w$ and taking $w\to\infty$ gives
\begin{equation}\label{e.A}
A(g_K)=- A(f^*) =  \frac 1{2\pi} \int_I \Im \bigl(f_I(x)\bigr) \,dx\,.
\end{equation}
Moreover,
$$
f^*(w)-w-2A(f^*)w^{-1}=
\frac 1{\pi}
\int_I\Im\bigl({f_I(x)}\bigr)\Bigl( \frac{1}{x-w}+\frac1w\Bigr)\, dx
\,,
$$
and therefore
\begin{eqnarray*}
\lefteqn {\bigl|f^*(w)-w-2A(f^*)w^{-1}\bigr|
} \\
&\le &
 \frac 1\pi \int_I 
\Im\bigl({f_I(x)}\bigr)\sup
\bigl\{|(x-w)^{-1}+w^{-1}|\st x\in I\bigr\}\, dx \\
&=&
 -2A(f^*) \sup\{ |x/((x-w)w)|\st x\in I\bigr\}\,.
\end{eqnarray*}
Hence, the proof will be complete once we demonstrate
that there is some constant $c_0$ such
that $I\subset [-c_0 r, c_0 r]$. This is easily done, as follows. 
 Define $G(z)\defeq g_K(rz)/r$ for $|z|>1$, and write 
    $G(z)=z+a_1 z^{-1}+ a_2 z^{-2}+\dots$.
The Area Theorem (see, e.g., \cite{Ru})
gives $1\ge \sum_{j=1}^\infty j |a_j|^2$.
In particular, $|a_j|\le 1$ for $j\ge 1$.
Consequently, we have $|G(z)-z|\le 1$ for $|z|\ge 2$.
By Rouch\'e's theorem (e.g.~\cite{Ru}), it follows that
$G\bigl(\{|z|\ge 2\}\bigr)\supset\{|z|>3\}$.  Consequently,
$g_K(\HH\setminus K)\supset\{|z|>3r\}$,
which gives $I\subset [-3r,3r]$.
\qed\bigskip

For convenience, we adopt the following notation
\begin{equation*}
K_{t,u}\defeq g_{K_t}(K_{t+u}\setminus K_t)
\, .
\end{equation*}

\medskip\noindent
{\bf Proof of Theorem \ref{driven}}.
We start with 1 implies 2.
Let $R\defeq \sup\{|z|\st z\in K_a\}$, and let
$Q\defeq\{z\in\HH\st |z|>R+2\}$.
Let $t,\delta,\epsilon$ and $S$ be as in the statement of the theorem,
and let $s\in \p S$.  Suppose that $\epsilon<1$ and
$r\in [\epsilon,\sqrt{\epsilon}]$.  Then there is
an arc $\beta_r$ of the circle of radius $r$
about $s$ such that $\beta_r\subset \HH\setminus K_t$
and $K_t\cup \R\cup \beta_r$ separates $K_{t+\delta}\setminus K_t$
from $Q$.  It therefore follows that the extremal length
of the set of arcs in $\HH\setminus K_t$ which separate
$K_{t+\delta}\setminus K_t$ from $Q$ in $\HH\setminus K_t$
is at most $\text{const}/\log (1/ \epsilon)$. 
(For the definition and basic properties of extremal
length, see \cite{A,LV}. 
The terms extremal length and extremal distance have the 
same meaning.) 
extremal length is invariant under conformal maps,
it follows that the extremal length of the set of
paths in $\HH$ that separate
$K_{t,\delta}$ from $g_{K_t}(Q)$
is at most $\text{const}/\log(1/\epsilon)$.
Because the diameter of $g_{K_t}(\HH\setminus Q)$ is bounded
by some function of $R$
(this follows since $g_{K_t}$ has the hydrodynamic normalization),
we conclude that at least one of these arcs has length less than
$\text{const}/\log(1/\epsilon)$.  Consequently,
this is a bound on the diameter of $K_{t,\delta}$.
Observe that this bound is uniform for all 
$t \in [0, a - \delta]$.
For each $t <a$, we  then define $U_t$ to be the point in the intersection
$\bigcap_{u>0}\overline {K_{t,u}}$.
We have an upper bound on 
$\diam( K_{t,\delta})$
which tends to zero uniformly as $\delta\to0$,
and therefore $\lim_{\delta\to0}g_{K_{t,\delta}}(z)-z=0$
uniformly for 
 $z\in\HH\setminus(K_{t,\delta})$ and 
$t \le a - \delta$. This implies that $U_t$ is uniformly continuous
on $[0, a )$ and can be extended continuously to $[0, a]$.
  
Now let $z_0\in \HH\setminus K_a$.
Then there is some $c>0$ such that $\Im\bigr(g_{K_t}(z_0)\bigl)>c$
for all $t\in[0,a]$.  Lemma \ref{cauchy} applied
with $K=K_{t,u}$, $z=g_{K_{t+u}}(z_0)$ and $x_0=U_{t+u}$ gives
$$
\frac{g_{K_t}(z_0)-g_{K_{t+u}}(z_0)}{u}
+\frac {2}{ g_{K_{t+u}}(z_0)-U_{t+u}}
\to 0
$$
as $\delta \to 0$.
As $g_{K_t}(z_0)$ and $U_t$ are  continuous in $t$, we may therefore
conclude
that 
$$
\p_t g_{K_t}(z_0)=\frac {2} { g_{K_t}(z_0)-U_t}\,,
$$
which gives 2.

The proof that 2 implies 1 is easy.
Let $\epsilon>0$.  Given $0\le t\le t+u < a$,
let $\rho(t,u)\defeq u + \max\{U(t')-U(t'')\st t',t''\in [t,t+u]\}$. 
Observe that $\diam(K_{t,u})\to 0$ if $\rho(t,u)\to 0$
and $\rho(t,u)\to 0$ if $u \to 0$.
Consequently, the extremal length of the set of paths
separating $K_{t,u}$ from $\{z\in\HH\st |z|>1\}$
in $\HH$ goes to zero as $\rho(t,u)\to0$. 
This implies that there is
a path $\beta$ in this set such that
$\diam\bigl(g_{K_t}^{-1}(\beta)\bigr)<\epsilon$, 
provided $u$ is small.
We then just take $S = g_{K_t}^{-1} (\beta)$.
\qed

\def\Timemodif/{Time-modified}
\def\timemodif/{time-modified}
\subsubsection {\Timemodif/ expanding hulls and  
restriction}

Let $(K_t,t\in[0,a])$ denote a family of hulls,
and suppose that there is a monotone increasing homeomorphism
$\phi:[0,a]\to[0,\hat a]$, such that
$\bigl(K_{\phi(t)},t\in[0,\hat a]\bigr)$ is an expanding hull
driven by some function $t\mapsto \hat U_t$.
If additionally $\phi$ is continuously differentiable in $[0,a]$
and $\phi'(t)>0$ for each $t\in[0,a]$,
then we call $(K_t,t\in[0,a])$ a \term{\timemodif/}
expanding hull, with driving function $U_t\defeq \hat U_{\phi^{-1}(t)}$.
Note that, in this case, $\phi^{-1}(t)=A(K_t)$, 
and that
\begin{equation}\label{e.smoothhull}
\p _t g_{K_t}(z) =\frac {2 \p _t A(K_t)}{g_{K_t}(z)-U_{t}}\,.
\end{equation}

Note that in our terminology, an expanding hull is always a \timemodif/ 
expanding hull.

\begin{lemma}\label{maphull}
Let $(K_t, t\in[0,a])$, be a \timemodif/ expanding hull,
with driving function $(U_t, t \in [0,a])$.
Let $D$ be a relatively open subset of $\closure\HH$ which contains
$\overline {K_a}$, and set $D_\R\defeq D\cap\R$. 
Let $G:D\to\closure\HH$
be conformal in $D\setminus D_\R$ and continuous in $D$,
and suppose that $G(D_\R)\subset\R$. 
Then $(G(K_t), t \in [0,a])$ is a \timemodif/ expanding hull.
Moreover, 
\begin{equation}\label{Ascal}
\p _t A\bigl(G(K_t)\bigr) =
G'(U_0)^2 \p _t A(K_t)\,,\qquad\text{at }t=0.
\end{equation}
\end{lemma}

\proof
We first prove (\ref{Ascal}).
The proof will be based on (\ref{e.A}).
Note first that if $K'=aK$ then $A(K')=a^2A(K)$.
Therefore, we may assume that $G'(U_0)=1$.
Similarly, with no loss of generality, we assume that
$U_0=G(U_0)=0$.
By the reflection principle, $G$ is analytic in $D$.

Set $\hat K_t\defeq G(K_t)$.
Let $I_t\subset \R$ be the interval 
corresponding to 
of $\closure{\p K_t\cap\HH }$ under $g_{K_t}$,
and let $\hat I_t$ be the
interval corresponding to 
$\closure{\p\hat K_t\cap\HH}$
under $g_{\hat K_t}$.  Let $\epsilon>0$, and
let $D_\epsilon\defeq \{z\in D\st |1-G'(z)|<\epsilon\}$.
Let $\beta$ be some arc in $D_\epsilon\setminus \{0\}$ that separates
$0$ from $\infty$ in $\closure \HH$.
Consider the map 
$$
h_t=g_{\hat K_t}\circ G \circ g_{K_t}^{-1}
.$$
It is well-defined in a neighborhood of $I_t$ 
provided that $K_t \cap \beta = \emptyset$ (for instance), and this holds
when $t$ is small. 
This map may be continued analytically by reflecting
in the real axis, and therefore the maximum
principle implies that 
$$
\sup\{h_t'(x)\st x\in I_t\}
\le
\sup \{|h_t'(z)|  \st
z\in g_{K_t}(\beta)\}
$$
when $K_t\cap\beta=\emptyset$.
Note that $g_{K_t}(z)-z\to 0$ and $g_{\hat K_t}(z)-z\to 0$
as $t\searrow 0$, and
therefore $g'_{K_t}(z)\to 1$ and $g'_{\hat K_t}(z)\to 1$
on $\beta$.  Consequently,  for small $t$
we have 
\begin{equation}\label{e.he}
\sup \{h_t'(x)\st x\in I_t\}< 1+2\epsilon.
\end{equation}
Note that for $z$ close to $U_0=0$ we have 
$\Im \bigl(G(z)\bigr)\le (1+\eps) \Im (z)$. 
Using (\ref{e.A}), this inequality and (\ref{e.he}), we get
\begin{eqnarray*}
A\bigl(G(K_t)\bigr)
&=&
 \frac 1{2\pi} 
\int _{\hat I_t} \Im\bigl( g_{\hat K_t}^{-1}(x)\bigr) \,dx
\\
&=&  \frac 1{2\pi} 
\int _{I_t} \Im\bigl( G\circ g_{K_t}^{-1}(x)\bigr) h_t'(x) \,dx
\\
& \le &  \frac 1{2\pi}\int_{I_t}
(1+\epsilon)
\Im\bigl(g_{K_t}^{-1}(x)\bigr) (1+2\epsilon) \,dx
\\
&=&
(1+\epsilon)
(1+2\epsilon) A(K_t)
\end{eqnarray*}
for small $t>0$.
(Note that $g_{K_t}^{-1}(x)$ is not defined for every $x\in {I_t}$, 
but it is defined for almost every 
 $x\in {I_t}$.)
By symmetry, we also have a similar inequality
in the other direction.  This proves (\ref{Ascal}).

By Theorem~\ref{driven}, to show that $G(K_t)$ is a \timemodif/
expanding hull, it suffices to show that
$A\bigl(G(K_t)\bigr)$ is continuously differentiable in $t$,
with derivative bounded away from $0$.
Let $G_t\defeq g_{\hat K_t}\circ G\circ g_{K_t}^{-1}$.
Then $G_t$ is analytic in $g_{K_t}(D\setminus K_t)$
and depends continuously on $t$.  Hence $G_t'(U_t)$ is
continuous in $t$.  Since
$A\bigl(G(K_{t+u})\bigr)=
A\bigl(G(K_t)\bigr) + A\bigl(G_t(K_{t,u})\bigr)$,
 it follows that
$\p _t A\bigl(G(K_t)\bigr) = G'_t(U_t)^2 \p _t A(K_t)$,
which completes the proof.
\qed\bigskip

For future reference, we note that when $g_t=g_{K_t}$ satisfies
the differential equation (\ref{e.smoothhull}), we have
the following formula 
\begin{equation}\label{eqder}
\p _t\log g_{t}'(z) = -\frac {2 \p _t
A(g_t)}{\bigl(g_{t}(z)-U_{t}\bigr)^2}\,,
\end{equation}
which is obtained by differentiating (\ref{e.smoothhull}) with respect
to $z$.

\subsubsection {Pairs of \timemodif/ expanding hulls}

We now discuss the situation where there are two disjoint expanding hulls.

Let  $(L_s, s\in [0,s_0])$, and $(K_t, t\in[0,t_0])$,
be a pair of \timemodif/ expanding hulls such that 
$\closure L_{s_0}\cap\closure {K}_{t_0}=\emptyset$.
Let $g_{s,t}\defeq g_{L_s\cup K_t}$,
$g_t \defeq g_{K_t}$, $\hat g_s\defeq g_{ L_s}$ and
$a(s,t)\defeq A(g_{s,t})$.
Then for each $s\in[0,s_0]$ and $t\in[0,t_0]$ we have
$$
g_{s,t}=g_{g_t(L_s)}\circ g_t
\,.
$$
Therefore, 
$$
\p _s g_{s,t}(z) =
 \frac{ 2\, \p _s a(s,t)}{g_{s,t}(z)-\Ub{s,t}}\,,
$$
where $s\mapsto \Ub{s,t}$ is the driving function for
the \timemodif/ expanding hulls $s\mapsto g_t(L_s)$.
Similarly,
$$
\p _t g_{s,t}(z) =
 \frac{ 2\, \p _t a(s,t)}{g_{s,t}(z)-\Ua{s,t}}\,,
$$
where $t\mapsto \Ua{s,t}$ is the driving function for
the \timemodif/ expanding hulls $t\mapsto \hat  g_s(K_t)$.
Although we do not know that $g_t^{-1}(\Ua{0,t})$ is well
defined, $g_{s,t}\circ g_t^{-1}$ is analytic in
a neighborhood of $\Ua{0,t}$, by the reflection principle.
Hence, it is clear that 
$\Ua{s,t} =  g_{s,t}\circ g_t^{-1}(\Ua{0,t})$ (see, for example,
the construction of $U_t$ in the proof of Theorem~\ref{driven}),
and therefore 
\begin{equation}\label{Uder1}
\p _s \Ua{s,t} = 
 \frac{ 2\, \p _s a(s,t)}{\Ua{s,t}-\Ub{s,t}}\,.
\end{equation}

We will now prove the formula
\begin{equation}\label{e.ader}
\p _s \p _t a(s,t) = 
 \frac {-4\, \p _ s a(s,t)\,\p _t a(s,t)}{(\Ua{s,t}-\Ub{s,t})^2}\,.
\end{equation}
{}From (\ref{Ascal}) we have 
$$
\p _t a(0,s)=  \hat g_s'(\Ua{0,0})^2 \p_t a(0,0) = \hat g_s'(\Ua{0,0})^2
\,,
$$
 and using~(\ref{eqder}), leads to 
\begin{eqnarray*}
{\p _s \log \p _t a(0,s)}
&=&
 \frac {-4 \,\p _s A(\hat g_s)}{\bigl(\hat g_{s}(\Ua{0,0})-\Ub{s,0}\bigr)^2}
\\
&=&
 \frac {-4 \,\p _s a(0,s)}{(\Ua{s,0}-\Ub{s,0})^2}\,.
\end{eqnarray*}
This verifies (\ref{e.ader}) for the case $t=0$.  The general case
is similarly obtained.

\subsection {Proof of Lemma \ref{l.loc}}\label{plocproofsec}

We will now prove Lemma \ref {l.loc}; this is the core
of the proof of the locality property. We that $\gamma : [0, s_1] 
\to \overline \HH$ is a continuous simple path 
with $\gamma (0) \in \R \setminus \{0\}$ and 
$L_s\defeq\gamma (0, s_1] \subset \HH$.
With no loss of generality, 
 assume that
$\gamma$ is parameterized so that
$A(L_s)=s$.
By Theorem~\ref{driven}, $(L_s , s \in [0, s_1])$ is a \timemodif/
expanding hull, and by Lemma~\ref{maphull},
for each $s$, $t\mapsto g_{L_s}(K_t)$
is a \timemodif/ expanding hull and
for each $t$, $s\mapsto g_{K_t}(L_s)$
is a \timemodif/ expanding hull.
Let $t\mapsto W(s,t)$ be the process driving
$t\mapsto g_{L_s}(K_t)$, let $U(s,t)$
be the process driving $s\mapsto g_{K_t}(L_s)$,
and let $Y_t=W(0,t)$ be the process driving $K_t$.
As above, let $a(s,t)=A(g_{s,t})$.
For simplicity, $W(s,t)$ will be abbreviated to
$W$, $U(s,t)$ to $U$, $a(s,t)$ to $a$, etc.

Our aim is to show that $(W(s_1, t)
, t \ge 0)$ is a continuous martingale (up to the stopping time $T$)
and that its quadratic variation (for background on stochastic calculus,
see, for instance, \cite {IW, RY}) 
is
$$
\langle W(s_1, \cdot ) \rangle_t
=
6 (a (s_1, t) - a (s_1, 0)).
$$
Indeed, if this is true,
let $\phi(t)$ be the inverse of the map
 $t\mapsto a(s_1, t )-a(s_1, 0)$, and
 define $ \tilde W (t) = W( s_1, \phi (t))$,
then $\tilde W (t/ 6)$ is a 
Brownian motion, so that 
$ t \mapsto g_{s_1, 0} (K_{\phi (t)}) - g_{s_1,0} (0)
$ 
is an \Ss/ process, as required.
Note that this will in fact give a precise expression for the
time-change in Lemma~\ref {l.loc} and Theorem~\ref {t.loc}.

\medbreak

Before giving the mathematically rigorous proof,
we first present a formal, nonrigorous
derivation of the fact that 
$W(s, \cdot)$ is a martingale. 
In this derivation,  
 $\slepar$ will
be kept as a variable, in order to stress where
the assumption $\slepar=6$ plays a role
(it will not be so apparent in our proof).

\medbreak
\noindent 
{\bf Nonrigorous Argument.}
The first goal is to show that the quadratic variation $\langle W\rangle_t$
of $t\mapsto W(s,t)$ satisfies
\begin{equation}
\partial_t \langle W\rangle_t = \slepar \partial_t a,   \label{e.spd}
\end{equation}
for each $s,t$. 
It is clear that this holds
when $s=0$, since $K_t$ is \Ss/.
We have
\begin{align*}
\partial_s\partial_t \langle W\rangle_t
&
= \partial_t\partial_s \langle W\rangle_t
\\&
= 2 \partial_t \langle\partial_s W,W\rangle_t
\\&
= 2 \partial_t
\langle2(\partial_s a)\bigl(W-U\bigr)^{-1},W\rangle_t
\qquad\qquad\hbox{(by (\ref{Uder1}))}
\\&
= -4 (\partial_s a)(W-U)^{-2}\partial_t \langle W\rangle_t
\qquad\qquad\hbox{(by It\^o's formula)}
\\&
= (\partial_t a)^{-1} (\partial_s\partial_t a)\partial_t \langle W\rangle_t
\qquad\qquad\qquad\hbox{(by (\ref{e.ader})).}
\end{align*}
Consequently,
$$
(\partial_t a)^2\partial_s\bigl(\partial_t \langle W\rangle_t  /\partial_t
a\bigr)
=\partial_t a \,\partial_s\partial_t \langle W\rangle_t -
\partial_t \langle W\rangle_t\,\partial_s\partial_t a
=0\,,
$$
which means that 
$\partial_t \langle W\rangle_t/\partial_t a$ does not depend on $s$.
Since (\ref{e.spd}) holds when $s=0$, 
 this proves (\ref{e.spd}).

We now show that $t\mapsto W(s,t)$ is a martingale.
The $dt$ term in It\^o's formula for the $\partial_t$-derivative of
$$
\partial_s W(s,t) = 2\partial_s a/(W-U)
$$
is
$$
\frac{2\partial_t\partial_s a}{W-U}
+2\frac{\partial_s a}{(W-U)^3} \partial_t\langle W\rangle_t
{}-4\frac{\partial_s a}{(W-U)^3}\partial_t a\,,
$$
where the  first summand comes from differentiating $\p _s a$,
the second summand is the diffusion term in It\^o's formula,
and the last summand comes from differentiating with respect to $U$ and
using
(\ref{Uder1}) for $\p _t U$.  Using (\ref{e.ader}) and (\ref{e.spd}),
this becomes
\[   \Bigl(3 - \frac{1}{2} \slepar\Bigr) \frac{\partial_t \partial_s 
  a}{W-U} \,  , \]
which vanishes when 
$\slepar=6$. Hence $t\mapsto\partial_s W(s,t)$ is a martingale.
As $W(s,t)=Y_t + \int_0^s \partial_s W(s',t)\,ds'$,
it follows that $t\mapsto W(s,t)$ is a martingale.
This completes the informal
 proof.

\medskip

The problem with the above argument is that we do
not know that $t\mapsto W(s,t)$ is a semi-martingale,
and hence cannot apply stochastic calculus to it.
Moreover, we need to check that there is sufficient
regularity to justify the equality
$\p _s \p _t \langle W\rangle_t=
\p _t \p _s \langle W\rangle_t$.

To rectify the situation, set
$$
V(s,t')\defeq
W(s,0)+\int_0^{t'} \sqrt{\p _t a(s,t)}\, d Y_t
\,.
$$
Then $t\mapsto V(s,t)$ is clearly a martingale.
The rest of this subsection will be devoted to the 
proof of the fact that $V=W$.
Recall that 
$$ T = \sup \{ t \ge 0 \st \overline K_t \cap \overline L_{s_1} = \emptyset
\} .$$
We will  need
the following fact:
\medbreak

\begin {lemma}
\label {Vcontinuity}
There exists a continuous 
version of $V$  on $[0, s_1] \times [0,T)$.
\end {lemma}

\proof
In order to keep some quantities bounded, we have to  
stop the processes slightly before
$T$.
Let us fix $\epsilon\in(0,1)$, and 
define
$$ 
T_1^\eps 
= \inf \{ t > 0 \ : \ \inf_{s \le s_1} 
| W(s,t) - U (s,t) | \le \eps \} .
$$

Define for any $s \le s_1$ and $t_0 \ge 0$.
$$ \tilde V (s, t_0) 
= \tilde V^{\eps} (s, t_0) 
\defeq \int_0^{\min (t_0, T_1^\eps)} \sqrt { \p_t a }  \,dY_t .
$$
As $\sup_{n \ge 1} T_1^{1/n} = T$, it is sufficient to
show existence of a continuous version (on $[0,s_1] \times \R_+$)
of $\tilde V$.

Let 
$$ \tilde a = \tilde a (s,t) 
 \defeq 1_{t \le T_1^\eps} \sqrt {\p_t a }.$$
 Note that $\p _t a(0,t)= \p _s a(s,0)=1$. 
 Hence, from (\ref{e.ader}) it follows that $\p _s a\le 1$ and $\p _t a\le
1$
 for all $s\le s_1$, $t<T$.
Using (\ref{e.ader}) again, we get
$$| \p_s   \tilde a |
=
1_{t \le T^\eps_1}
 \frac { 2 \sqrt {\p_t a } \p_s a } { (W-U)^2 } 
\le
2  \eps^{-2} .$$
Hence,  
for all $t \ge 0$, for all $s,s'$ in $[0, s_1]$,
$$
| \tilde a (s, t) - \tilde a (s', t) | \le {2  \eps^{-2}}
 |s-s'|.
$$
But
\begin {eqnarray*}
\lefteqn {
\expect \left[ (\tilde V(s, t_0) - \tilde V(s', t_0'))^4 \right]}
\\
&\le &
16\, \expect \left[ \Bigl( \int_{t_0}^{t_0'} \tilde a (s', t) \,dY_t 
\Bigr)^4 \right]
+
16\, 
\expect \left[\Bigl
( \int_0^{t_0} \bigl(\tilde a (s,t) - \tilde a (s', t)\bigr) \,dY_t \Bigr)^4
\right] 
\end {eqnarray*}
and using, for instance, the 
Burkholder-Davis-Gundy inequality for $p=4$
(see, e.g., \cite[IV.4]{RY}),
 we see that 
there exists a constant
$c_1=c_1(\eps)$ such that for all $t_0, t_0' \ge 0$
and  $s,s' \in [0,s_1]$ 
\begin {eqnarray*}
\lefteqn {
\expect \left[ \Bigl(\tilde V(s, t_0) - \tilde V(s', t_0')\Bigr)^4 \right]}
\\
&\le &
c_1 
\expect \bigl[(t_0 - t_0')^2\bigr] +
c_1 \expect \left[ \Bigl(\int_0^{t_0} (s-s')^2 \,dt \Bigr)^2 \right] 
\\
& \le &
c_1  (t_0- t_0')^2 + c_1 t_0^2 (s-s')^4 
\end {eqnarray*}
and the existence of a continuous version of $\tilde V$ then easily
follows
from Kolmogorov's lemma (see, e.g., \cite[I.(1.8)]{RY}).
\qed

\medskip

{}From now on, we will use a version of $V$ that is continuous on
$[0,s_1] \times [0, T)$.
Define 
\begin {eqnarray*}
T_2^\eps &\defeq& \inf \{ t \ge 0 \ : \ 
\sup_{s \le s_1} | V (s,t) - W(s,t) | \ge 1 \}\,, \\
T_3^\eps &\defeq& \inf \{ t \ge 0 \ : \ 
\inf_{s \le s_1} | V(s,t) - U(s,t) | \le  \eps \}\,, 
\\
T^\eps &\defeq& \min (T_1^\eps, T_2^\eps, T_3^\eps )\,
. \end {eqnarray*}

Note that
for all  $s \le s_1$ and $t < T$,
$$
 \partial_s   \sqrt {\partial_t a }
=
 \frac { -2 \sqrt {\partial_t a} \partial_s a }{(W - U)^2 } 
.$$
The process $\partial_s \sqrt {\partial_t a } $
remains bounded before $T^\eps$ (uniformly in $s \le s_1$),
and  it
is a measurable function of $(s,t)$. 
By Fubini's Theorem for
stochastic integrals 
 (see \cite[Lemma III.4.1]{IW}),
we have that for all $s_0\le s_1$,
for all $t_0 \ge 0$, 
 almost surely
\begin {eqnarray}
\lefteqn {\int_0^{s_0} \int_0^{t_0'}
 \frac { -2 \sqrt {\partial_t a} \partial_s a }{(W - U)^2 } 
 \,dY_t\, ds}  \nonumber 
\\
&=&
\int_0^{t_0'}  \Bigl(\int_0^{s_0} \p_s  \sqrt{\p_t a}\,  ds \Bigr)\,dY_t 
\nonumber \\
&=&
\int_0^{t_0'} 
\left( 
\sqrt {\p_t a (s_0, t) } - \sqrt {\p_t a (0, t) } 
\right) \,dY_t 
\nonumber
\\
&=&
V(s_0, t_0') - V(0, t_0' ) -\bigl( W(s_0, 0) - W(0,0)\bigr) 
\label {eq1}
\end {eqnarray}
where $t_0' \defeq \min (t_0, T^\eps)$.
On the other hand, using It\^o's formula, we now compute 
\begin{equation}
\label {eq2}
\begin {aligned}
{ \frac  {2\partial_s a (s,t_0')  }{  U(s,t_0') - V(s,t_0') }
}
\ =\ &
\frac {2} {U(s, 0) - V(s, 0) }
+ \int_0^{t_0'}
\frac {2\p_s a \sqrt{\p_t a}}{(U-V)^2}\, dY_t
\\ &\quad
{} + 2
\int_0^{t_0'}
\left(
\frac { \p_t \p_s a }{U-V} - \frac {\p_s a\, \p_t U }{(U-V)^2 } 
+ \frac {\p_s a\, \p_t \langle V \rangle_t }{(U-V)^3}
\right) 
\,dt
\\
\ =\ & - \partial_s W(s,0)
+
\int_0^{t_0'}
{\tilde b_1} \,dY_t + \int_0^{t_0'} (V-W) b_2 \, dt 
\end {aligned}
\end{equation}
where (using $\p _t \langle V\rangle_t = \p_t a\,\p_t \langle Y\rangle_t = 
\slepar\, \p_t a=6\,\p_ta$)
\begin {eqnarray*}
\tilde b_1 (s,t) &\defeq& 
\frac{2 \p_s a \sqrt {\p_t a} }{(U-V)^2 }\\
b_2 (s,t) &\defeq&
\frac {4\p_s a\, \p_t a }{ (U-W)^2 (U-V)^2 }  
\left( 5+ 3 \frac {W-V}{V-U} \right)
.\end {eqnarray*}
Note that for all $s \le s_1$ and $t \le T^\eps$,
$$ | b_2 (s,t) | \le {16} \eps^{-5}\,. $$
By integrating (\ref {eq2}) with respect to $s$
and subtracting 
(\ref {eq1}) from it, we get 
\begin {eqnarray}
\lefteqn {V(s_0, t_0') - V(0, t_0')
 + \int_0^{s_0} 
\frac {2\partial_s a (s,t_0')  }{ U(s,t_0') - V(s,t_0') }
\,ds 
}
\nonumber
\\
&=&
\int_0^{s_0} 
\int_0^{t_0'} (V-W)  b_2 (s,t) \,dt\, ds 
\nonumber
\\
&& \ +
\int_0^{s_0} 
\left(  \int_0^{t_0'} \frac { - 2 \sqrt {\partial_t a} \partial_s a 
}{ (W - U)^2 }\, dY_t  + \int_0^{t_0} \tilde b_1 (s,t) \,dY_t \right)  \,ds
\nonumber
\\
&=&
\int_0^{s_0}
\int_0^{t_0'}
(V-W) b_2 \,dt \,ds + 
\int_0^{s_0} \int_0^{t_0'}
 (V-W) b_1  \,dY_t  \,ds 
\label {eq3}
\end {eqnarray}
where (after some simplifications)
$$ b_1 (s,t) 
\defeq
\frac {2\, \p_s a \sqrt {\p_t a}}{ (V-U)^2 
(W-U)^2 } \bigl((U-W)+(U-V)\bigr).
$$  
Note that for all $s \in [0, s_1]$ and $t \le T^\eps$,
$$ | b_1 (s,t) | \le 4\eps^{-3}  \,.
$$
But we know on the other hand that
\begin {equation}
\label {eq4}
W(s_0, t_0') - W(0, t_0')
 =  \int_0^{s_0}
\frac  {2\,\partial_s a (s,t_0') }{  W(s,t_0') - U(s,t_0')  }
\,ds .
\end {equation}
Subtracting this equation from (\ref {eq3}),  one gets
\begin {eqnarray*}
\lefteqn {V(s_0, t_0') - W(s_0, t_0')
 = \int_0^{s_0} b_3 (s,t_0') (V(s,t_0')-W (s,t_0')) \,ds
}\\
&&
+
\int_0^{s_0}
\int_0^{t_0'}
(V-W) b_2 \,dt \,ds + 
\int_0^{s_0} \int_0^{t_0'}
 (V-W) b_1  \,dY_t  \,ds 
\end {eqnarray*}
where
$$ b_3 (s,t) \defeq \frac {-2\, \p_s a }{ (U-W) ( U-V) }.
$$
Again $b_3$ remains uniformly bounded before $T^\eps$.

We now define
$$ H(s,t)= V(s,t) - W(s,t) .$$
Hence, for all $t_0 \le T$,
$$
H(s_0, t_0)
= \int_0^{s_0} 
b_3 (s, t_0) H(s, t_0) \,ds
+
 \int_0^{s_0} \int_0^{t_0} 
b_2 H \,ds \,dt
+ 
\int_0^{s_0} \int_0^{t_0} 
b_1 H \,ds \,dY_t 
$$
and $|b_1|, |b_2|, |b_3|$ are all bounded by some
constant $c_2=c_2(\eps)$ on $[0, s_1] \times [0, T^\eps]$.
This equation and
an argument similar to Gronwall's Lemma will
show that $H=0$.

Let us fix $t_1 > 0$.
For any $t \ge 0$, define 
$$ \tau (t) = \min (t, t_1, T^\eps).$$
We will use the notation $\tau_0 =  
\tau (t_0)$.
It is easy to see that there
exists a $c_3 = c_3 ( \eps, t_1, s_1)$ such that 
for all $t_0 \ge 0$ and $s_0 \in [0, s_1]$,
\begin{multline}
H(s_0, \tau_0)^2 
\le 
c_3  \int_0^{s_0} H(s, \tau_0)^2  \,ds \\
{}+ c_3 \int_0^{\tau_0} \int_0^{s_0} H^2 \,ds \,dt
+  c_3\Bigl(\int_0^{\tau_0}\int_0^{s_0} b_1 H \,ds \,dY_t \Bigr)^2 .
\end{multline}
The Burkholder-Davis-Gundy
inequality for $p=2$  (see, e.g., \cite[IV.4]{RY})
 shows that
there exists constants $c_4= c_4 (\eps, t_1, s_1)$
and $c_5 = c_5 (\eps, t_1, s_1)$ such that for all
$s \in [0,s_1]$ and $t_0 \ge 0$,
\begin {eqnarray*}
\expect
\left[ \sup_{u \le t_0 }
\Bigl(\int_0^{\tau (u)} \int_0^{s_0} b_1 H \,ds \,dY_t \Bigr)^2 \right]
&\le & c_4  \,
\expect \left[ \int_0^{\tau_0} \Bigl(\int_0^{s_0} b_1 H \,ds \Bigr)^2  \,dt
\right] 
\\
& \le &
c_5 \int_0^{t_0}  \int_0^{s_0}
\expect\left[ H \bigl(s, \tau (t)\bigr)^2\right] \,ds \,dt 
.\end {eqnarray*}

Let us now define
$$ h(s_0, t_0 ) = 
\expect \left[ \sup_{t \le t_0 } H\bigl( s_0, \tau (t) \bigr)^2  \right].$$

Then
$$
h(s_0, t_0) 
\le
c_6  \left( 
\int_0^{s_0} h(s, t_0) \,ds 
+ \int_0^{t_0} \int_0^{s_0} h (s,t)  \,ds \,dt 
\right).
$$
We also know that 
$h(s,t)$ is bounded by $1$
(because $|H| \le 1$). Hence, it is straightforward to prove by induction
that for all $s_0 \in [0, s_1]$,  $t_0 \ge 0$
and $p=1,2,\dots$,
$$ h(s_0, t_0) 
 \le \frac { c_6^p s_0^p (1+t_0)^p }{ p! } \,,
$$ 
so that $h (s_0, t_0) =0 $.
In particular (using the continuity of $V $ and $W$),
 this shows that 
$ W = V $
almost surely on all  sets 
$[0, s_1] \times [0, \min (t_1, T^\eps )]$.
As this is true for all $\eps$ and $t_1$, we 
conclude that 
$V=W$ on $[0, s_1] \times [0, T)$.
Lemma~\ref{l.loc} follows, and thereby also Theorem~\ref{t.loc}.
\qed

\section {Exponents for \Ss/}
\label{s.Sexp} 

\subsection {Statement}\label{setup}

In the present section, we are going to compute
intersection exponents associated with \Ss/.

Suppose that $D\subset\C$ is a Jordan domain; that is,
$\p D$ is a simple closed curve in $\C$.
Let $a,b\in \p D$ be two distinct points on the
boundary of $D$.  As explained in Section~\ref{slebasic},
the \Ss/ $(K_t, t \ge 0)$ from $a$ to $b$ in $D$ is well defined,
up to a linear time change.

Now suppose that $I\subset\p D$ is an arc with $b\in I$
but $a\notin I$.  Let 
$$
\tau_I\defeq\sup\{t\ge 0\st \closure {K_t}\cap I=\emptyset\}\,.
$$
By Corollary~\ref{split}, up to a time change,
the law of the process $(K_t,  t<\tau_I)$
does not change if we replace $b$ by another point
$b'\in I$.
Set 
$$
S=S({a,I,D})\defeq \bigcup_{t<\tau_I} K_t
\,,
$$
and call this set the \term{hull from $a$ to $I$ in $D$}.
It does not depend on $b$.

Suppose that $L >0$, and  let ${\cal R}={\cal R}(L)$ 
denote the rectangle with corners
\begin{equation}\label{corners}
A_1\defeq 0\,,\qquad 
A_2\defeq L\,,\qquad
A_3\defeq L+i\pi\,,\qquad
A_4\defeq i\pi\,.
\end{equation}
Let ${\cal S}$ denote the closure of
the hull from $A_4$ to $[A_1,A_2]\cup[A_2,A_3]$ in ${\cal R}$.

In the following, we will use the terminology
\term{$\pi$-extremal distance} 
instead of ``$\pi$ times the extremal distance''.
For instance, the $\pi$-extremal distance between the
vertical sides of ${\cal R}$ in ${\cal R}$  
is $L$.

When ${\cal S} \cap [A_1,A_2] = \emptyset$, let
${\cal L}$ be the $\pi$-extremal distance  
 between $[A_1, A_4]$
and $[A_2,A_3]$ in ${\cal R} \setminus {\cal S}$.
Otherwise, put ${\cal L} = \infty$.
 
In the sequel, we will use the function
\begin{equation}\label{udef}
u(\lambda) = \frac { 6 \lambda + 1
 + \sqrt { 24 \lambda +1}  } { 6} \,.
\end{equation}

The main goal of this section is to prove the following
result.  

\begin {theorem}
\label {expS}
\begin {equation}
\label {explambda}
\expect \left[  1_{ {\cal L} < \infty }
\exp ( - \lambda {\cal L} ) \right]
=\exp\Bigl( -  u(\lambda)  L + O(1)(\lambda+1)\Bigr),
\qquad\hbox {as } L \to \infty\,, 
\end {equation}
for any $\lambda \ge 0$ 
(where $O(1)$ denotes an arbitrary quantity whose absolute value is bounded
by a constant which does not depend on $L$ or $\lambda$).
\end {theorem}

\medbreak

In particular, when $\lambda=0 $,
\begin {equation}
\label {epx1/3}
\prob \bigl[ {\cal S} \cap [A_1, A_2] = \emptyset \bigr] 
 = \prob [ {\cal L} < \infty ] =
\exp \bigl( - L / 3  + O (1) \bigr) ,\qquad \hbox {as } L \to \infty.
\end {equation}

\subsection {Generalized Cardy's formula}

By conformal invariance, we may work in
the half plane $\HH$.
Map the rectangle ${\cal R}$ conformally
onto $\HH$, so that $A_1$ is mapped to $1$,
$A_2$ is mapped to $\infty$, $A_3$ is mapped to
$0$, and then the image $x=x(L)\in(0,1)$
of $A_4$ is determined for us.
Let $K_t$ be the hull of an \Ss/ process $g_t=g_{K_t}$ in
$\HH$, with driving process $W(t)$, which is started at $W(0)=x$.
(That is, $K_t$ is a translation by $x$ of the standard
\Ss/ starting at $0$.)
In order to emphasize the dependence on $x$, we will
use the notation $\prob_x$ and $\expect_x$ for 
probability and expectation.

Set
\begin{equation*}
\begin{aligned}
T_0 &\defeq \sup
\{ t \ge 0 \st \closure K_t \cap (-\infty,0] = \emptyset \}\,,
\\
T_1 &\defeq \sup
\{ t \ge 0 \st \closure K_t \cap [1,\infty) = \emptyset \}\,,
\\
T &\defeq \min\{T_0,T_1\}\,.
\end{aligned}
\end{equation*}
As will be demonstrated, $T_0,T_1<\infty$ a.s.
Let
$$
f_t(z) \defeq
\frac { g_t (z) - g_t (0) }{g_t (1) - g_t (0)}\,,
$$
for $t<T$, which is just $g_t$ renormalized to fix
$0,1$ and $\infty$.  It turns out that
$f_T\defeq \lim_{t\nearrow T} f_t$ exists a.s.
On the event $T_1>T_0$, $f_T$ uniformizes
the quadrilateral 
$$\bigl(\HH\setminus K_T; 1,\infty,\min(\closure K_T\cap\R),
\max(\closure K_T\cap \R)\bigr)
$$ to the
form
$$\Bigl(\HH; 1,\infty,0,f_T\bigl(\max(\closure K_T\cap \R)\bigr)\Bigr).$$
Therefore, we want to know the distribution of
$$1 - f_T ( \max ( \closure K_T \cap \R ))$$
and how it depends on $x$ (especially when 
$x \in (0,1)$ is close to 1).
In this subsection, we will calculate something very closely related:
the distribution of $f_T' (1)$ and how it depends on $x$.

Set
$$
\Lambda(1-x,b)\defeq
\expect_{x} 
\bigr[ {1}_{\{T_0<T_1\}}\, f'_T (1)^b \bigl]
$$
for $b\geq 0$ and $x \in (0,1)$.
Recall the definition of the hypergeometric function  
${}_2F_1$ (see, e.g., \cite {Le}):
$$
 {}_2F_{1} (a_0, a_1, a_2 ;\, x   ) =
 \sum_{n=0}^\infty \frac {(a_0)_n(a_1)_n}{(a_2)_n n!}\, x^n\,,
 $$
 where $(a)_n = \prod_{j=1}^n (a+j-1)$ and $(a)_0 = 1$.
Note that ${}_2F_1(a_0, a_1,  a_2 ;\, 0)=1$.

\begin {theorem} 
\label {cardy}
For all $b \ge 0$, $x \in (0,1)$, 
$$
\Lambda(x,b)
=
\frac {\sqrt{\pi}\,2^{-2\hat b}\Gamma(5/6+\hat b)
}{\Gamma(1/3)\Gamma(1+\hat b)}\,
x^{1/6+\hat b}\,
{}_2F_1(1/6+\hat b,1/2+\hat b,1+2\hat b;\,x)
$$
where
$$
\hat b=\frac {\sqrt{1+24b}}{6}.
$$
and ${}_2F_1$ is the hypergeometric function.
\end{theorem}

\medbreak

Setting $b=0$, we obtain Cardy's formula \cite{CaFormula},
as in \cite{S2}.  Thus, this 
result can be thought of as a
generalization of Cardy's formula.

\medbreak

Note that Theorem~\ref{cardy} determines  completely 
the law of 
$1_{T_0<T_1} f_T' (1)$.
In particular, the Laplace transform of 
the conditional law of $\log 1/ f_T' (1)$ 
given $\{ T_0<T_1 \}$
is 
$\Lambda (1-x,b) / \Lambda (1-x,0)$.

\proof
We first observe that $T_1<\infty$ almost surely.
Indeed, $g_t(1)-W(t)$ is a Bessel process with index $5/3$, with
time linearly scaled,
and hence hits $0$ almost surely in finite time (e.g.~\cite {RY}).
Similarly, $g_t(0)-W(t)$ hits $0$ almost surely. 
It is clear that $T_1$ is the time $t$ when $g_t(1)-W(t)$ hits
$0$ and $T_0$ is the time when $g_t(0)-W(t)$ hits $0$.
It follows from Theorem~\ref{driven} and $x\in(0,1)$
that almost surely, $T_0\neq T_1$. 
We may also  conclude that
\begin{equation}
\label{hitfirst}
\lim_{x\to0}\prob_x[T_0<T_1] =1\,,
\qquad
\lim_{x\to1}\prob_x[T_0<T_1] =0\,.
\end{equation}

The next goal is to prove that
\begin{equation}\label{deriv}
\lim_{t\nearrow T} f_t'(1)>0\quad \hbox{if and only if}\quad T_0<T_1\,.
\end{equation}
If $T_0<T_1$, then $\closure K_T\cap [1,\infty]=\emptyset$.
Therefore $f_T$ is defined and conformal near $1$,
and $f'_T(1)>0$, by the reflection principle.
On the other hand, if $T_1<T_0$, then
$\closure K_T\cap[1,\infty)\neq\emptyset$. 
We claim that 
\begin{equation}\label{surround}
1\notin\closure K_{T_1}\hbox{ almost surely;}
\end{equation}
since $\closure{K_{T_1}}\cap[1,\infty)\neq\emptyset$,
this means that $K_{T_1}$ separates $1$ from $\infty$
in $\HH$. 
Indeed, let $\phi:\HH\to\HH$ be the anti-conformal automorphism
that fixes $x$ and exchanges $1$ and $\infty$.
$T_1<\infty$ a.s.\ and $\sup_{t<T_1}|W(t)|<\infty$ a.s.\ imply that
$K_{T_1}$ is bounded a.s., which is the same as
saying that $\phi(K_t)$ stays bounded away from $1$
as $t\nearrow T_1$.
But Corollary~\ref{split} and invariance under reflection imply
that up to time $T_1$ the law of $K_t$ is the same as a time-change
of the law of $\phi(K_t)$. 
Hence, a.s.\ $K_t$ stays bounded away from $1$ as $t\nearrow T_1$,
proving (\ref{surround}).
It follows from (\ref{surround}) that $\lim_{t\nearrow T_1} f_t'(1)=0$
a.s.\ on the event $T_1<T_0$ (observe that,
given (\ref{surround}), the extremal length from a neighborhood
of $0$ to a neighborhood of $1$ in $\HH\setminus K_t$ tends to
$\infty$ as $t\nearrow T_1$), and (\ref{deriv}) is established.

Define the renormalized version of $W(t)$: 
$$
Z (t)\defeq \frac  { W(t) - g_t (0) }{g_t (1) - g_t (0)}  \,,
$$
and the new time-parameter
$$
\uu = \uu(t) 
\defeq \int_0^t \frac { dt }{ (g_t (1) - g_t (0) )^2 },\qquad t<T\,.
$$
Set $\uu_0\defeq \lim_{t\nearrow T}\uu(t)$.
Since $T_1\neq T_0$ a.s., $\inf\{g_t (1) - g_t (0)\st t<T\}>0$
a.s., and hence $\uu_0<\infty$ a.s.
Let $t(\uu)$ denote the inverse to the map $t\mapsto \uu(t)$.
A direct calculation gives
$$
\partial_{\uu} (f_{t(\uu)} (z))
=
\frac {-2}{Z(t)- f_t(z)} 
+ 
\frac {2 (1-f_t (z) ) }{Z(t)} - \frac {2 f_t (z) }{1- Z(t)}
\,,
$$ 
and
$$
dZ (t) =
\frac { dW_t }{ g_t(1) -g_t (0) } 
+
\frac {2 dt }{ (g_t(1) - g_t (0))^2 }
\left( \frac {1} {Z(t)} + \frac {1} {Z(t) - 1} \right). $$
We now use the notation
$$
\tilde Z (\uu) \defeq  Z( t(\uu)) ,\qquad 
\tilde f_\uu (z) \defeq f_{t(\uu)} (z)\,.
$$
Then,
\begin {equation}
\label {evolZ}
d\tilde Z_\uu = 
dX_\uu + 
\frac {2 (1- 2 \tilde Z (\uu)) \,d\uu}{ \tilde Z (\uu) 
(1- \tilde Z (\uu))}
=
dX_\uu + 
 \left( \frac {2}{ \tilde Z (\uu) } - \frac {2} { 1 -  \tilde Z (\uu) }
\right) d\uu \,,
\end {equation}
where $(X_\uu, \uu \ge 0)$ has the same law
as $(W(t), t \ge 0)$; i.e., it is a 
Brownian motion with time rescaled by a factor of $6$.
Also,
\begin {equation}
\label {evolf'}
\partial_\uu (\tilde f_\uu (z))
= \frac { -2 }{ \tilde Z (\uu) - \tilde f_\uu (z)}
+
\frac {2 (1-\tilde f_\uu (z) ) }{\tilde Z(\uu)}
 - \frac {2 \tilde f_\uu (z) }{1- \tilde Z(\uu)}.
\end {equation}
These two equations describe the evolution of $\tilde f_\uu (z)$.
Note that $\uu(T)=\uu_0$ is the first time at which 
$\tilde Z (\uu) $ hits 0 or 1.

We now assume that $b > 0$.
Differentiating (\ref {evolf'}) with respect to $z$ gives
(the Cauchy integral formula, for example, shows that we
may indeed differentiate, but this is also
legitimate since $\p_\uu$ and $\p_z$ commute
in this case)
\begin {equation}
\label {evolder}
\partial_\uu (\log \tilde f_\uu' (z))
=
\frac {-2 } { (\tilde Z(\uu) - \tilde f_\uu (z) )^2 }
   - \frac {2 }{ \tilde Z (\uu) } - 
\frac {2 } {1- \tilde Z (\uu) }.
\end {equation}
We are particularly interested in 
$$ \alpha (\uu) \defeq   \log \tilde f_\uu' (1) = \log f_{t(\uu)}' (1)$$
which satisfies
\begin {equation}
\label {evolalpha}
\partial_\uu \alpha (\uu)
=
\frac {-2 } { (\tilde Z(\uu) - 1 )^2 }
   - \frac {2 }{ \tilde Z (\uu) } -
\frac {2 } {1- \tilde Z (\uu) }.
\end {equation}
Note that equations (\ref {evolZ}) and (\ref {evolalpha})
describe the evolution of 
the Markov process $\bigl(\tilde Z(\uu), \alpha (\uu)\bigr)$.
The process stops at $\uu_0$.
Define
$$
y(x,v)
\defeq  \expect \left[ \exp (b \alpha (\uu_0)) \mid \tilde Z (0)=x , \alpha
(0)
= v \right]
$$
where the expectation corresponds to the Markov process 
started from $\tilde Z (0)=x$ and $\alpha (0)=v$.
{}From the definition of $y$ it follows that
$$
y(x,0)=\Lambda(1-x,b)\,,
$$
since $\lim_{t\nearrow T} f_t= f_T$ in a neighborhood of $1$
on the event $T_0<T_1$ and (\ref{deriv}) holds.
It is standard that such a function $y(x,v)$ is $C^\infty$, and
the strong Markov property
ensures that the process 
$$Y= y\bigl( \tilde Z (\uu), \alpha (\uu) \bigr)$$
is a local martingale.
The drift term in It\^o's formula for $d Y$ must vanish,
which gives
\begin {equation}
\label {PDE}
0
=
\frac {2 (1-2x) }{x (1-x)} \p_x y
+ 3 \p_{xx}^2 y 
+ 
\left( 
\frac {-2 } {(1-x)^2 }
   - \frac {2} {x}
   - \frac {2} {1-x}
\right)
\p_{v} y .
\end {equation}
As 
$$ 
\alpha (\uu) = \alpha (0) 
+ \int_0^u 
( \partial_\uu \alpha (\uu'))
\,d\uu' ,$$
we get
$$ y(x,v) = \exp (bv) y(x,0).$$
Set
$$ h ( x ) \defeq y (1-x , 0)=\Lambda(x,b),$$
so that $y (x,v) = \exp (bv) h(1-x)$. Hence (\ref {PDE})
becomes
\begin {equation}
\label {ODE}
-2b h(x) + 2x (1-2x) h'(x) + 3 x^2 (1-x) h'' (x) = 0
.\end {equation}
The second statement in (\ref{hitfirst})
implies that
$$
\lim_{x\searrow 0} h(x)=0\,,
$$
while
$$
\lim_{x\nearrow 1} h(x)=1
$$
holds, since when $x$ is close to $0$,
$K_{T_0}$ is likely to be small, by scale invariance, for example.
The differential equation (\ref{ODE})
 can be solved explicitly by looking for solutions
of the type
$
h(x) = x^{c} z(x)$:
two linearly independent 
solutions are ($i=1,2$) 
$$
h_i (x) =
x^{1/6 + b_i} 
 {}_2F_1(1/6+ b_i,1/2+ b_i ,1+2b _i ;\, x ) 
$$
where
$$ 
b_1 = - b_2 = \frac {\sqrt {1+ 24 b} } {6} 
.$$
Recall that ${}_2F_1(a_0, a_1,  a_2 ;\, 0)=1$.
The function $h(x)$ must be a linear 
combination of $h_1$ and $h_2$.
However, $\lim_{x\searrow 0}h(x)=0=\lim_{x\searrow 0} h_1(x)$, 
but $\lim_{x \searrow 0} h_2 (x)=\infty $ .
Hence, $h(x)=c h_1(x)$ for some constant $c$.
The equality $h(1) = 1$ and knowledge of the
value at $x=1$ of hypergeometric functions
(see, e.g., \cite {Le}) allows the determination of $c$, and establishes
the theorem in the case $b>0$.  The case $b=0$ follows
by taking a limit as $b\searrow 0$.
\qed

\medskip

\noindent{\bf Remark.}
With the same proof, Theorem~\ref{cardy} generalizes to \SG/ with
$\slepar>4$,
and gives 
$$
\Lambda_{\slepar}(x,b)
=
C(b,\slepar)
x^{1/2-2/\slepar+\hat b_\slepar}
{}_2F_1(1/2-2/\slepar+\hat b_\slepar,
    6/\slepar-1/2+\hat b_\slepar,1+2\hat b_\slepar;\,x)\,,
$$
where
$$
\hat b_\slepar\defeq \frac {\sqrt{(\slepar-4)^2 + 16 \slepar b}}{2
\slepar}\,,
\qquad 
C(b,\slepar)\defeq
\frac
{ \Gamma(3/2-6/\slepar+\hat b_\slepar) \Gamma(1/2+2/\slepar+\hat b_\slepar)}
{\Gamma(1-4/\slepar)\Gamma(1+2\hat b_\slepar)}\,.
$$

\subsection {Determination of the \Ss/ exponents}

For every $t\geq 0$, set 
$$
M_t\defeq \max(\closure K_t\cap \R)\,.
$$
The following Lemma shows that 
our understanding of the derivative $f_T'(1)$ gives information
on $f_T(M_T)$ itself. 

\begin{lemma}
\label {theta}
In the above setting, let $N_T\defeq  f_T(M_T)$. 
For $b\geq 0$, set
$$
\Theta(x,b)\defeq 
E_{1-x}\bigl[{1}_{\{T_0<T_1\}}\,(1-N_T)^b\bigr].
$$
Then
\begin{equation}
\label{sandwich}
\left({x/ 2}\right)^b
\Lambda(x/2,b) \leq \Theta(x,b)\leq x^b \Lambda(x,b).
\end{equation}
\end{lemma} 

Note that $\Theta(x,b)$ is close 
to the quantity we are after,
since $-\log(1-N_T)$ is approximately the extremal 
length  of the quadrilateral
$\bigl(\HH\setminus K_T;
\min(\closure K_T\cap\R),\max(\closure K_T\cap\R),1,\infty\bigr)$. 

\proof
It follows easily from (\ref{evolder}) that 
$f'_t(z)$ is nondecreasing in $z$ (viewed
as a real variable),
 as long as $z\geq M_t$.
Therefore,
$$
1-N_T=
\int_{M_T}^1  f_T'(z)\,dz
\leq (1-M_T) f_T'(1)
\leq x f_T'(1).
$$
This gives the right hand inequality
in (\ref{sandwich}).

To get the other inequality, consider some fixed $x^*>1$
($x^*$ should be thought of as 
close to 1; we will eventually take
$x^* = 1 / (1-x)$).
Let $\tilde x^*= f_T(x^*)$.
Then
$$
\tilde x^*- N_T\geq
\tilde x^*-1=\int_{1}^{x^*}  f_T'(z)\,dz
\geq (x^*-1)  f_T'(1).
$$
This gives
\begin{equation}\label{e.xo}
E_{1-x} \bigl[{1}_{\{T_0<T_1\}}(\tilde x^*-N_T)^b\bigr]
\geq (x^*-1)^b\Lambda(x,b).
\end{equation}
A simple scaling argument will give
an upper bound of the left-hand side of this inequality in terms
of $\Theta $.
Let 
$$T^*
= \inf \{ t \ge 0 \st
\closure K_t\cap\R\setminus(0 ,x^*) \neq\emptyset \} . $$
Note that $T \leq T^*<\infty$ a.s.\ and that 
$T=T^*$ if $T_0<T_1$.
For each $t\leq T^*$, let $f^*_t$ be the conformal map from the
unbounded component of $\HH \setminus K_t$ to $\HH$,
 which fixes the
points $\infty$, $0$, $x^*$. For all $t\le T$,
$$f^*_t (z) =\frac {x^* }{f_t(x^*)} f_t (z) .$$
Note that 
$$f_t (x^*) \le x^* = f_t^* ( x^*) .$$
Then,
$$
{1}_{\{T_0<T_1\}}(\tilde x^*-N_T)
\leq 
{1}_{\{T_0<T_1\}}\bigl(x^*-f_T^*(M_T)\bigr).
$$
Hence,
\begin{equation}\label{e.sll}
{1}_{\{T_0<T_1\}}(\tilde x^*-N_T)
\leq 
{1}_{\{M_{T^*}< x^*\}}\bigl(x^*-f_{T^*}^*(M_{T^*})\bigr),
\end{equation}
since on the event $T_0<T_1$, we have $T=T^*$ and
$M_T=M_{T^*}$.
However, by scale invariance,
when $W(0)=1-x$ the random variable
$$ {1}_{\{M_{T^*}< x^*\}}\bigl(x^*-f_{T^*}^*(M_{T^*})\bigr) $$
has the same law as the random variable
$$x^* {1}_{\{M_T< 1\}}\bigl(1- f_T(M_T)\bigr) $$
does when $W(0)=(1-x)/x^*$.
Thus, combining (\ref{e.sll}) and (\ref{e.xo}) gives
\begin{equation*}
(x^*)^b\Theta\bigl(1-(1-x)/x^*,b\bigr)\geq(x^*-1)^b\Lambda(x,b).
\end{equation*}
We take $x^*=1/(1-x)$, say, and get
$$\Theta(2x,b)\geq\Theta(2x-x^2,b)\geq x^b\Lambda(x,b),
$$
which gives the left hand side of (\ref{sandwich}).
\qed

\medbreak

\noindent
{\bf Proof of Theorem \ref {expS}.}
We are now ready to derive Theorem \ref {expS}
by combining Theorem \ref {cardy}
and Lemma \ref {theta}.

In the setting of the theorem, let $\phi:{\cal R}(L)\to\HH$ be the
conformal homeomorphism satisfying $\phi(A_1)=1,\phi(A_2)=\infty$,
and $\phi(A_3)=0$.  Set $x=x(L)\defeq\phi(A_4)$.
By conformal invariance, the law of ${\cal L}$
is the same as that of the $\pi$-extremal distance ${\cal L}^*$
from $(-\infty,0]$ to $(N_T,1)$ in $\HH$ (with the notations
of Lemma~\ref{theta}). Considering the map
$z\mapsto \log (z-1)$ makes it clear that
\begin{equation}\label{L}
L= - \log(1-x) +O(1)
\,,
\end{equation}
for $L>1$, and similarly
\begin{equation}
\label{Lstar}
{\cal L}^* = -\log (1-N_T) +O(1)
\,.
\end{equation}
For $L >1$, (note also that ${\cal L} \ge L$),
\begin{equation*}\begin{aligned}
\expect \left[  1_{ {\cal L} < \infty }
\exp ( - \lambda {\cal L} ) \right]
&
=
\expect \left[  1_{ {\cal L}^* < \infty }
\exp ( - \lambda {\cal L}^* ) \right]
\\&
=
\exp\bigl(O(1)\bigr)
\expect \left[  1_{N_T<1} (1-N_T)^{ \lambda } \right]
\qquad\hbox{(by (\ref{Lstar}))}
\\&
=
\exp\bigl(O(1)\bigr)
\Theta(1-x,\lambda)
\\&
=
\exp\bigl(O(\lambda+1)\bigr)
(1-x)^{u(\lambda)}
 \qquad\hbox{(by \ref{theta} and \ref{cardy})} 
\\&
=
\exp\bigl(O(\lambda+1)\bigr)
\exp\bigl( - {u(\lambda)L}\bigr)
 \qquad\hbox{(by (\ref{L}))} 
\,,
\end{aligned}\end{equation*}
which completes the proof of
Theorem~\ref {expS}. \qed

\section {The Brownian half-plane exponents}
\label{s.Bexp} 

We are now ready to combine the results collected 
so far and a ``universality'' idea 
similar to that developped in \cite {LW2} to 
compute the exact value 
of some Brownian intersection exponents 
in the half-plane.

\subsection {Definitions and background}

In this short subsection, we 
quickly review some results on 
intersection exponents between independent planar Brownian motions.
For details and complete proofs of these results, see \cite {LW1, LW2}.

Suppose that $n+p$
independent planar Brownian motions
$\beta^1, \cdots , \beta^n$ and $
\gamma^1, \cdots , \gamma^p$ are started from 
points $\beta^1(0) = \cdots = \beta^n(0) = 0  $ and 
$\gamma^1(0) = \cdots = \gamma^p(0) = 1$ in the complex plane,
 and 
consider 
the probability $f_{n,p}(t)$ that for all $j \le n$ and $l\le p$,
the paths of $\beta^j$ up to time $t$
and of
$\gamma^l$ up to time $t$  do not intersect;
more precisely:
$$ f_{n,p} (t)\defeq 
\prob \left[\Bigl( \bigcup_{j=1}^n \beta^j [0,t]\Bigr) \cap
\Bigl(\bigcup_{l=1}^p \gamma^l [0,t] \Bigr) = \emptyset \right].$$
It is  easy to see that as $t \to \infty$ this probability decays 
roughly 
like a power of $t$. The $(n,p)$-intersection exponent
$\xi (n,p)$ is defined as twice this power, i.e.,
$$ f_{n,p} (t) = (\sqrt {t})^{- \xi (n,p)  +o(1) },  \qquad t \to
\infty\,.$$
We call $\xi (n,p)$  the intersection exponent between
one packet of $n$ Brownian motions and one packet of
$p$ Brownian motions (for a  list of references 
on Brownian intersection exponents, see \cite {LW1}).
Note that the exponent $\zeta$ described in the introduction is 
$\xi (1,1)/2$.
It turns out to be more convenient to use this definition 
as a power of $\sqrt {t}$, i.e., of the space parameter. 
A Brownian motions travels very roughly to distance $\sqrt {t}$ in time $t$:
recall that if $\beta$ is a planar Brownian motion started from
$0$, say, and $T_R$ denotes its hitting time of the
circle of radius $R$ about $0$, then 
for all $\delta >0$, the probability that 
$T_R \notin (R^{2-\delta}, R^{2+ \delta})$
decays as $R\to\infty$ faster than any negative power of $R$.
This facilitates an easy conversion between the time based definition of
intersection exponents and a definition where the particles
die when they exit a large ball.

Similarly, one can define corresponding probabilities
for intersection exponents in a half-plane
$$ \tilde f_{n,p} (t)\defeq 
\prob \left[
\forall j\le n ,\,\forall l\le p,\,\,\,
\beta^j[0,t]\cap\gamma^l[0,t]=\emptyset\hbox{ and }
\beta^j[0,t]\cup\gamma^l[0,t]\subset {\cal H}\right],
$$
where ${\cal H}$ is some half-plane containing the two starting points.
($\tilde f_{n,p}(t)$ will depend on ${\cal H}$.)
In plain words, we are looking at the probability that all
Brownian motions stay in the half-plane and that all $\beta$'s 
avoid all $\gamma$'s.
It is also easy to see that there exists a $\tilde \xi (n,p)$ 
(which does not depend on ${\cal H}$)
such that 
$$
\tilde f_{n,p} (t) = (\sqrt {t})^{-\tilde \xi (n,p) +o(1)},\qquad
t\to\infty\,.
$$
Note that  $\tilde \zeta$ described
in the introduction is $\tilde \xi (1,1) /2$.

One can also define intersection exponents
$\xi (n_1, \ldots , n_p)$ and $\tilde \xi (n_1, \ldots, n_p)$
involving more 
packets of Brownian motions.
(For a more detailed discussion of this see \cite {LW1}).
For instance, if $B^1, B^2, B^3,  B^4$
denote four Brownian motions started from different points, 
the exponent $\xi (2, 1, 1)$ is defined 
by
\begin{multline*}
\prob
\left[
\hbox {The three sets }
B^1[0,t] \cup B^2 [0,t] 
, \ 
B^3[0,t] 
, \
B^4[0,t] 
\hbox { are disjoint}
\right] 
\\
=
t^{- \xi (2,1,1) / 2 +o(1) },\qquad t\to\infty\,.
\end{multline*}
One of the results of \cite {LW1} is that there is a natural
and rigorous way to generalize the definition of 
intersection exponents
between packets of Brownian motions to the case where each packet
of Brownian motions
is the union of a ``non-integer number'' of paths;
for the half-plane exponents, one can define the exponents 
$\tilde \xi ( u_1, \ldots, u_p )$, where
$u_1, \ldots, u_p  \ge 0$.
These generalized exponents satisfy the so-called cascade
relations (see \cite {LW1}):
for any 
$1 \le q \le p-1$,
\begin {equation}
\label {cascade}
\tilde  \xi ( u_0, \ldots, u_p)
=
\tilde   \xi (u_0, \ldots, u_{q-1}, \tilde \xi (u_{q},
\ldots , u_p))
. \end {equation}
Moreover, $\tilde\xi$ is invariant under a permutation
of its arguments.

\medbreak

There exists (see \cite {LW1, LW2})  a characterization
 of these
exponents in terms of the 
so-called Brownian excursions 
that turns out to be 
useful.
For any bounded simply connected open domain $D$, there exists a
Brownian excursion measure $\mu_D$ in $D$.
This is an infinite measure on paths $(B(t) , t \le \tau)$
in $D$ such that $B(0,\tau) \subset D$ and
$B(0), B(\tau) \in \partial D$ (these can viewed as prime
ends if necessary).
$x_s\defeq B(0)$ and 
$x_e\defeq B(\tau)$ are the starting point and terminal point of the
excursion.
One possible definition of $\mu_D$
 is the following: Suppose first that 
$D$ is the unit disc. For any $s>0$ define the 
measure $P^s$ on Brownian paths (modulo continuous increasing time-change)
 started uniformly on the 
circle of radius $\exp (-s)$, and killed when they exit $D$.
Note for any $s_0> s$,
  the killed 
Brownian path defined under the probability measure
$P^s$  has a probability 
$s/ s_0$ to intersect the circle of radius $\exp (- s_0)$.
Then, define $$\mu_D \defeq 
\lim_{s \searrow 0} (2 \pi  r /  s)  P^s .
$$
One can then easily check 
that for any M\"obius transformation $\phi$ from $D$ onto $D$, 
  $ \phi ( \mu_D) = \mu_D$.
This makes it possible to extend the definition of $\mu_D$
to any 
simply connected domain $D$, by conformal invariance.
These Brownian excursions  
also have a ``restriction''  property 
\cite {LW2}, as the Brownian paths only feel the
boundary of $D$ when they hit it (and get killed).

Suppose for a moment that ${\cal R}={\cal R}(L)\subset\C$
is the 
rectangle with corners given by (\ref{corners}), and that ${\cal B}$ is the
trace
of the Brownian excursion $(B(t), t \le \tau)$ in ${\cal R}$.
Define the event
$$
E_1 = 
\{ B(0) \in [A_1,A_4] \hbox { and } B(\tau) \in [A_2,A_3] \}\,,
$$
i.e., $B$ crosses the rectangle from the left to the right.
(Although $\mu_{{\cal R}}$ is an infinite measure,
$\mu_{{\cal R}}(E_1)$ is finite.)
When $E_1$ holds, let ${\cal R}_B^+$ be the
component of ${\cal R}\setminus {\cal B}$ above ${\cal B}$,
and let ${\cal R}_B^-$ be the 
component of ${\cal R}\setminus {\cal B}$ below ${\cal B}$.
Let ${\cal L}_B^-$
(respectively ${\cal L}_B^+$)
denote the $\pi$-extremal distance between $[A_1, x_s] $ and $[A_2, x_e ]$
in ${\cal R}_B^-$
(respectively $[x_s, A_4 ] $ and $[x_e, A_3]$)
in ${\cal R}_B^+$.

Then, for any $\alpha \ge 0$ and $\alpha' \ge 0$,
the exponent $\tx (\alpha, 1, \alpha') = 
\tx (1, \tx (\alpha, \alpha'))$ is characterized by 
\begin{equation}\label{honey}
\expect_{\mu_{\cal R}}
\left[
1_{E_1}
\exp ( - \alpha {\cal L}^+_B - \alpha' {\cal L}^-_B ) 
\right]
=
\exp ( - \tx (\alpha', 1, \alpha ) L + o(L) ) \,,
\end{equation}
when $L \to \infty$, where $ \expect_{\mu_{\cal R}}$
denotes expectation (that is, integration)
with respect to the measure ${\mu_{\cal R}}$.
Similarly,
\begin{equation}\label{plop}
\expect_{\mu_{\cal R}}
\left[
1_{E_1}
\exp ( - \alpha {\cal L}^+_B ) 
\right]
=
\exp ( - \tx (1, \alpha ) L + o(L) ) \,,\qquad L\to\infty\,.
\end{equation}
See~\cite{LW2}.
It will also be important later that $\tx$ is continuous in
its arguments, and that $\lambda\mapsto \tx(1,\lambda)$ is
strictly monotone.

\subsection {Statement and proof}

For any $p \ge 0$, we put 
$$ v_p = \frac {p (p+1) }{6}
.$$
Let
${\cal V}$ denote the set of numbers
$ \{  v_p  \st 
p \in \N \} 
$.
Note that the smallest values in ${\cal V}$ are
$0, 1/3, 1, 2, 10/3, 5, 7$.

We are now ready to prove the following result:

\begin {theorem}
\label {general}
For any $ k \ge 2$,
$\alpha_1, \cdots , \alpha_{k-1}$
in ${\cal V}$ 
 and for all $\alpha_k \in \R_+$, 
\begin {equation}
\label {eqgeneral}
\tx (\alpha_1, \cdots , \alpha_k )
=
\frac {\left( 
\sqrt { 24 \alpha_1 + 1} +
 \cdots + \sqrt { 24 \alpha_k +1 } 
   - (k-1) \right)^2  - 1  }{ 24 } .
\end {equation}
\end {theorem} 
 
\medskip

It is immediate to verify that this Theorem implies
Theorem~\ref{theoremintro}.
\medskip

\noindent {\bf Remark.}
In \cite{LSW3} Theorem~\ref{general} is extended to all nonnegative
reals $\alpha_1, \ldots, \alpha_k$.

\medskip

Theorem \ref {general} is a 
consequence of the cascade relations and the following lemma, which is the
special case of the theorem with $k=2, \alpha_1
=1/3$:

\begin {lemma}
\label {uu}
For any $\lambda > 0$,
$$
\tx ( 1/3  , \lambda ) 
=
u ( \lambda)\,,
$$
where $u(\lambda)$ is given by {\rm (\ref{udef})}. 
\end {lemma}

\noindent {\bf Proof of Theorem \ref  {general}}
(assuming Lemma \ref {uu}).
Define for all $\lambda \ge 0$, 
$ U(\lambda) = \sqrt {24 \lambda + 1 } - 1 $.
Lemma \ref{uu} implies immediately that for all $\lambda \ge 0$, 
$$
U( \tx (1/3, \lambda ) )
= U(\lambda ) + 2 
=U ( \lambda ) +  U(1/3) 
$$
and (for all integer $p$), 
$v_{p+1} = \tx (1/3, v_p)$.
The cascade relations then imply that
 for all 
integers  $p_1, \ldots , p_{k-1}$, 
$$ 
\tx ( v_{p_1} , \ldots , v_{p_{k-1}} , \lambda )
= 
U^{-1} \big(  2 (p_1 + \ldots + p_{k-1} )  + U (\lambda)\big). 
$$
This is  
(\ref {eqgeneral}).
\qed

\medbreak
\noindent
{\bf Proof of Lemma \ref {uu}.} 
For convenience, we again  work in a rectangle rather than in the
upper half-plane.
Let ${\cal R} = {\cal R} (L)$, 
and let ${\cal S}$ denote the closure of the hull of \Ss/ from
$A_4$ to $[A_1,A_2]\cup[A_2,A_3]$ in ${\cal R}$, as in~\ref{setup}.
Let $ {\cal B} $ denote the trace of a Brownian excursion 
in ${\cal R}$; we will
call its starting point $x_s$ and its terminal point $x_e$.
Consider the following events:
\begin {eqnarray*}
E_1 & = & \{ x_s \in [A_1,A_4 ] 
\hbox { and } x_e \in  [ A_2,A_3 ]  \}\,, \\
E_2 & = & \{ {\cal S } \cap [A_1,A_2 ] = \emptyset \}\,,
\\
E_3
&=& E_1 \cap E_2 \cap \{ {\cal S} \cap {\cal B} = \emptyset \}
\,.\end {eqnarray*}
When $E_2$ holds, let ${\cal L}_S$ denote the $\pi$-extremal
distance between the vertical edges of ${\cal R}$ in ${\cal R} \setminus
{\cal S}$
(that is, in the quadrilateral \lq\lq below\rq\rq\ ${\cal S}$).
Otherwise, let ${\cal L}_S = \infty$.

When $E_1$ holds, let ${\cal R}_B^+$ be the
component of ${\cal R}\setminus {\cal B}$ above ${\cal B}$,
and let ${\cal R}_B^-$ be the 
component of ${\cal R}\setminus {\cal B}$ below ${\cal B}$.
Let ${\cal L}_B^-$ (respectively ${\cal L}_B^+$)
denote the $\pi$-extremal distance between the vertical edges of ${\cal R}$
in ${\cal R}_B^-$ (respectively in ${\cal R}_B^+$), as before.
When $E_3$ holds, let 
${\cal L}_{SB} $ denote the  $\pi$-extremal
distance between the vertical edges of ${\cal R}$
in ${\cal R }_B^+ \setminus {\cal S}  $
(that is, in the quadrilateral \lq\lq below ${\cal S}$ and above
${\cal B}$\rq\rq).

\medbreak

Let $\lambda> 0 $.
We are interested in 
the asymptotic behavior of 
$$
f(L) = \expect \big[ 1_{  E_3 } \exp ( - \lambda {\cal L}_{SB} ) \bigr] 
$$
when $L \to \infty$.
By first taking expectations with respect to $B$ (with
the measure $\mu_{\cal R}$), 
and using the restriction property (Cor.~\ref{rest})
 for the domains
${\cal R}$ and ${\cal R}_B^+$,
it follows that as $L\to\infty$,
\begin{eqnarray*}
\label{chain1}
f(L)
&=
& \expect_{B}
\Bigl[ \expect_{S} \bigl[ \exp ( - \lambda {\cal L}_{SB} )\bigr]\Bigr]
\\&
=
&\expect_B \Bigl[ \exp \bigl( - u(\lambda) 
{\cal L}_B^+ +O(1) \bigr) \Bigr]
\qquad\hbox{(by Thm.\ \ref{expS} and restriction to ${\cal R}_B^+$)}
\\&
=
& \exp \Bigl( -  \tx \bigl(1, u( \lambda )\bigr) L +o(L) \Bigr)
\qquad\hbox{(by (\ref{plop}))}
 \,.
\end{eqnarray*}
On the other hand, we may first take
expectation with respect to ${\cal S}$.  Given ${\cal S}$,
the law of ${\cal L}_{SB}$ is the same as that
of ${\cal L}_B^-$, by complete conformal invariance of
the excursion measure (which is the analogue of the
restriction property to the excursion measure; see \cite{LW1}). 
Hence, as $L\to\infty$,
\begin{eqnarray*}
f(L) &=
& \expect_S \left[ \expect_B \left[  
 1_{E_3} \exp ( - \lambda {\cal L}_{SB} )
 \right]\right]
\\ &
= 
& \expect_S \left[ \expect_B \left[  
 1_{E_3} \exp ( - \lambda {\cal L}_{B}^- )
 \right]\right]
\\ &
= 
& \expect_B \left[ \expect_S \left[  
 1_{E_3} \exp ( - \lambda {\cal L}_{B}^- )
 \right]\right]
\\
&=
& 
 \expect_B \left[ \prob_S[E_3\mid {\cal L}_B^+]
 \exp ( - \lambda {\cal L}_{B}^- )
 \right]
\\ &
=
& \expect_B \left[ \exp\bigl(-{\cal L}_B^+/3+O(1)\bigr)
 \exp ( - \lambda {\cal L}_{B}^- )
 \right]
\qquad\hbox{(by (\ref{epx1/3}))}
\\ &
=
&\exp\Bigl(-\tx\bigl(1/3,1,\lambda\bigr)L + o(L)\Bigr)
\qquad\hbox{(by (\ref{honey}))}
\\ &
=
&\exp\Bigl(-\tx\bigl(1,\tx(1/3,\lambda)\bigr)L +o(L)\Bigr)\,,
\end{eqnarray*}
by the cascade relations (\ref{cascade}).
Comparing with (\ref{chain1}) gives
$\tx (1, \tx (1/3, \lambda) ) = \tx (1, u(\lambda)) $.
Finally,
$$
\tx (1/3,  \lambda ) = u( \lambda)\,,
$$
follows, since $\lambda'\mapsto \tx(1,\lambda')$ is strictly
increasing.  
\qed

\section{Crossing exponents for critical percolation}
\label{s.percolation} 

It has been conjectured \cite{S1} that \Ss/ corresponds
to the scaling limit of critical percolation clusters.
As additional support for this conjecture, we now show that
it implies the conjectured formula
for the exponents corresponding to the probability
that a long rectangle is crossed by $p$ disjoint
paths or clusters 
of critical percolation (\cite {DS2, Ca2, ADA}).

Let us first explain the conjectured relation between
\Ss/ and critical percolation.
Let $D\subset\C$ be a domain whose boundary $\p D\subset\C$
is a simple closed curve.
Let $a,b\in\p D$ be distinct points. Let
$\gamma_1$ be the counterclockwise arc on $\p D$ from $a$ to $b$,
and let $\gamma_2$ be the clockwise arc on $\p D$ from $a$ to $b$.
Let $\delta>0$, and consider a fine hexagonal grid
$H$ in the plane with mesh $\delta$; that is,
each face of the grid is a regular hexagon with
edges of length $\delta$, and each vertex has
degree $3$.  For simplicity, assume that $\p D$
does not pass through a vertex of $H$ and
that $a$ and $b$ do not lie on edges of $H$.
Color each hexagon of $H$ independently,
black or white, with probability $1/2$.
Then the union of the black hexagons forms one of the standard
 models for critical percolation
(see Grimmett~\cite{Gr} for percolation background and
references).

Given the random coloring, there is a unique path $\beta\subset\closure{D}$
that starts at $a$, ends at $b$, such that 
whenever $\beta$ is not on $\gamma_1$ it has
a black hexagon on its ``right'' and whenever $\beta$
is not on $\gamma_2$ it has a white hexagon on
its ``left''. 
This path is the boundary between the union of the
white clusters in $D$ touching $\gamma_2$
and the black clusters in $D$ touching $\gamma_1$. 
  Let $f: D\to\HH$ be a conformal homeomorphism
such that $f(a)=0$ and $f(b)=\infty$,
and parameterize $\beta$ in such a way that $A\bigl(f(\beta[0,t])\bigr)=t$.
Let $D_t$ be the component of $D\setminus \beta[0,t]$ that has $b$
on its boundary, and let $K_t=D\setminus D_t$.
The conjecture from \cite{S1}
(stated a bit differently) is that as $\delta\to 0$ the process
$(K_t, t\ge 0)$ converges to \Ss/ from $a$ to $b$ in $D$.
In light of this conjecture, the Locality Theorem \ref{t.loc}
and its corollaries are very natural.

Now consider an arc $I\subset \p D$, which contains $b$ but not $a$. 
Let $b_1$ and $b_2$ be the endpoints of $I$, labeled in such a way
that the triplet $a,b_1,b_2$ is in counterclockwise order around $D$.
Let $\gamma_1'\subset\gamma_1$ be the counterclockwise arc from $a$ to
$b_1$,
and let $\gamma_2'\subset\gamma_2$ be the clockwise arc from $a$ to $b_2$.
Let $T$ be the first time such that $\beta(t)\in I$, and set
$S\defeq \bigcup_{t<T} K_t$.
Then the component $\alpha_1$ of $\p S\cap D$ joining $\gamma_1'$ to $I$ is
a crossing in $\closure{{\cal B}}$ from $\gamma_1'$ to $I$,
which is \lq\lq maximal\rq\rq, in the sense that any other crossing
$\alpha\subset \closure{{\cal B}}$ from $\gamma_1'$
to $I$ is separated by $\alpha_1$ from $b_2$ in $D$.

Let $L$ be large, and
recall the definition of the rectangle  ${\cal R}={\cal R}(L)$
with corners given by~(\ref{corners}).
Let $p\in\N_+$, and
$\sigma=(\sigma_1,\dots,\sigma_p)\in\{\hbox{black},\hbox{white}\}^p$.
Consider the event $C_\sigma({\cal R})$ that there are paths
$\alpha_1,\dots,\alpha_p$, from $[A_4,A_1]$ to $[A_2,A_3]$ in
${\cal R}$ such that each $\alpha_j$ is contained in the union of
the hexagons of color $\sigma_j$, there is no hexagon which intersects
more than one of these paths, and $\alpha_{j+1}$ separates $\alpha_{j}$
from $[A_1,A_2]$ in ${\cal R}$ when $j=1,2,\dots,p-1$.

Take $\alpha_1$ to
be the topmost crossing with color $\sigma_1$, if such
exists, let $\alpha_2$ be the topmost crossing with color
$\sigma_2$ which is below all the hexagons meeting
$\alpha_1$, etc.  Then $C_\sigma({\cal R})$ holds iff these specific
$\alpha_1,\dots,\alpha_p$
exist.
Note that after we condition on $\alpha_1$, the hexagons
``below'' it are still independent and are black or white with
probability $1/2$.
Hence the following formula holds:
$$
\prob[C_{\sigma}({\cal R})]=
 \prob[\alpha_1\hbox{ exists}]\,
 \expect\Bigl[\prob\bigl[ C_{\sigma_{+1}}({\cal R}_{\alpha_1})
 \mid\alpha_1\bigr]\,\Big|\,
 \alpha_1\hbox{ exists}\Bigr]\,,
$$
where $\sigma_{+1}=(\sigma_2,\sigma_3,\dots,\sigma_p)$ ,
${\cal R}_{\alpha_1}$ is the union of the hexagons below $\alpha_1$,
and $C_{\sigma_{+1}}({\cal R}_{\alpha_1})$ is the event
that there are multiple crossings with colors specified by
$\sigma_{+1}$ from $[A_4,A_1]$ to $[A_2,A_3]$ in
${\cal R}_{\alpha_1}$.
 
It is clear that $\prob[C_\sigma({\cal R})]$ does not depend on the
choice of the sequence $\sigma$, but only its length.
Moreover, the conjectured conformal invariance
(or the conjecture that \Ss/ is the scaling limit)
implies that $\lim_{\delta\to 0} \prob[C_\sigma({\cal D})]$,
depends on the quadrilateral ${\cal D}$ only through its conformal modulus.
Hence define
$$
f_p(L)\defeq 
\lim_{\delta\to 0} \prob\Bigl[C_\sigma\bigl({\cal R}(L)\bigr)\Bigr]\,,
\qquad \sigma\in\{\hbox{black},\hbox{white}\}^p\,.
$$
We also set $f_p(\infty)\defeq 0$ and $f_0(L)\defeq 1_{L<\infty}$.

Let $S$ be the \Ss/ hull from $A_4$ to $I\defeq [A_1,A_2]\cup [A_2,A_3]$
in ${\cal R}={\cal R}(L)$, as defined in Subsection~\ref{setup}.
Let ${\cal R}_-$ be the component of ${\cal R}\setminus S$
which has $A_1$ on its boundary, and let
${\cal L}$ denote the $\pi$-extremal length from $[A_4,A_1]$
to $[A_2,A_3]$ in ${\cal R}_-$.
Note that ${\cal L}=\infty$ if $\closure S \cap [A_1,A_2]\neq\emptyset$.
Then we have
$$
f_p(L)=\expect [f_{p-1}({\cal L})],\qquad p=1,2,\dots
\,.
$$
To completely justify this step requires more work,
which we omit, since this whole discussion depends
on a conjecture anyway.  The slight difficulty
has to do with the fact that having a crossing 
of a closed rectangle is a closed 
condition, and the probability of a closed event
can go up when taking a weak limit of measures.
One simple way to deal with this is to note that
when the continuous process has a crossing in the rectangle
${\cal R}(L+\eps)$, every sufficiently close discrete
approximation of it has a crossing of the rectangle
${\cal R}(L)$.

Consequently, induction and Theorem~\ref{expS} give
\begin {equation}
\label {expperc}
f_p(L)= \exp\bigl(- (L+O(1))v_p\bigr)\,, \qquad L\to\infty
\,,
\end {equation}
where
\begin {equation}
\label {vp}
v_p = u^{\circ p} ( 0) = \frac { p (p+1)}{6}\,,
\end {equation}
 as before.
Here, the constant implicit in the $O(1)$ notation may depend on $p$.

Note also that if 
$$\sigma = ( \hbox {white, black, white , black,} \ldots , \hbox {white})
\in \{\hbox{black},\hbox{white}\}^{2k-1} \,,$$
then (in the discrete setting), the 
event $C_\sigma ({\cal R})$ is identical to the event
that the rectangle ${\cal R}$ is crossed from
left to right by $k$ disjoint white clusters.

The exponents 
(\ref {vp}) are  
those  predicted in \cite{Ca2,ADA}.

\vskip 1cm

G.L. :

Department of Mathematics

Box 90320

Duke University

Durham NC 27708-0320, USA

jose@math.duke.edu

\medbreak

O.S. : 

Microsoft Research

1, Microsoft Way

Redomond WA 98052, USA 

schramm@microsoft.com

\medbreak

W.W. :

D\'epartement de Math\'ematiques

B\^at. 425

Universit\'e Paris-Sud

91405 ORSAY cedex, France

wendelin.werner@math.u-psud.fr

\end {document}